\documentclass{article}

\usepackage{arxiv}

\usepackage[utf8]{inputenc} % allow utf-8 input
\usepackage[T1]{fontenc}    % use 8-bit T1 fonts
\usepackage{hyperref}       % hyperlinks
\usepackage{url}            % simple URL typesetting
\usepackage{booktabs}       % professional-quality tables
\usepackage{amsfonts}       % blackboard math symbols
\usepackage{nicefrac}       % compact symbols for 1/2, etc.
\usepackage{microtype}      % microtypography
\usepackage{lipsum}
\usepackage{graphicx}
\graphicspath{ {./images/} }

\usepackage{amssymb,amsmath,amsthm}

\usepackage{tikz-cd}

\newcommand{\R}{\mathbb{R}}
\newcommand{\C}{\mathbb{C}}
\newcommand{\CP}{\mathbb{P}}
\newcommand{\Z}{\mathbb{Z}}

\DeclareMathOperator{\re}{Re}
\DeclareMathOperator{\im}{Im}

\newcommand{\mc}{\mathcal}

\DeclareMathOperator{\id}{id}

\title{Böttcher-type potential for the secant map}

\author{
 Nicholas Freeman \\
  Department of Mathematics\\
  University of Michigan\\
  Ann Arbor, MI 48109 \\
  \texttt{njfree@umich.edu} \\
}

\begin{document}
\maketitle
\begin{abstract}
We present a construction of a Böttcher-type holomorphic map for the potential of the secant method dynamical system near a root-type fixed point. The modulus of the Böttcher-type map extends to be continuous on the entire basin of attraction of the fixed point, and is real-analytic away from the iterated preimages of the fixed point. Using this construction, we show the associated Green's function for the fixed point is pluriharmonic wherever it is finite.
\end{abstract}

\section{Introduction}
The secant method is a root-finding algorithm for single-variable smooth or holomorphic functions. In this article, we take a connected open set $U\subseteq\C$ and consider a nonconstant meromorphic function $f : U\dashrightarrow \C$. Then the secant method to solve $f(z) = 0$ for $z$ takes the form
\[
z_{n+1} = z_n - \frac{f(z_n)(z_n-z_{n-1})}{f(z_n)-f(z_{n-1})} = \frac{f(z_n)z_{n-1} - f(z_{n-1})z_n}{f(z_n)-f(z_{n-1})}.
\]
This iteration is typically presented as an alternative to Newton's method. It is well-known that for any pair $(z_0,z_1) \in U\times U$ of initial guesses that are sufficiently close to a simple root $\alpha\in U$ of $f$, the associated sequence $(z_n)_{n\in\Z_{\geq0}}$ converges to $\alpha$ superlinearly. More precisely,
\[
\lim_{n\to\infty} z_n = \alpha \text{ implies } \lim_{n\to\infty} \frac{|z_{n+1}-\alpha|}{|z_n-\alpha|^\phi} \leq \lambda
\]
for some constant $\lambda\geq 0$ independent of $z_0$ and $z_1$, where $\phi = \frac12(1+\sqrt5)$ is the golden ratio. See Díez \cite{DiezSecantConvergence} for a development of this convergence fact.

One can treat the secant method iteration as a fixed-point iteration of the secant map $\mc{S}_f : U\times U\dashrightarrow \CP^1\to \CP^1$ given by
\[
\mc{S}_f(x,y) = \left(\frac{f(x)y-f(y)x}{f(x)-f(y)}, x\right)
\]
in the standard coordinates of $\CP^1\times \CP^1$, so that $\mc{S}_f(z_n,z_{n-1}) = (z_{n+1},z_n)$ for all $n\in\Z_{\geq0}$. This $\mc{S}_f$ is holomorphic and sensibly extends to equal the graph of Newton's method at most points along the diagonal of $U\times U$. In particular, Bedford and Frigge \cite{BedfordFriggeSecant} show that if $z_0\in U$ is a simple root of $f$, then the diagonal point $(z_0,z_0)$ is a superattracting fixed point of $\mc{S}_f$ with complex Jacobian
\[
D\mc{S}_f(z_0,z_0) = \begin{bmatrix}
    0 & 0\\
    1 & 0
\end{bmatrix}.
\]
In particular, standard facts of dynamical systems theory then implies there is a basin of attraction
\[
\mc{A}(f,z_0) = \left\{(x,y)\in U\times U : \lim_{n\to\infty} \mc{S}_f^{\circ n}(x,y) = (z_0,z_0)\right\}
\]
which is nonempty, open, and fully invariant under $\mc{S}_f$, i.e., we have $\mc{S}_f(x,y)\in \mc{A}(f,z_0)$ if and only if $(x,y)\in \mc{A}(f,z_0)$.

Following Hubbard and Papadopol, \cite{HubbardPapadopolPotential} we seek to understand the potential function $h_{f,z_0} : \mc{A}(f,z_0) \to \R$ given by
\[
h_{f,z_0}(x,y) = \lim_{n\to\infty} \left\|\mc{S}_f^{\circ n}(x,y) - (z_0,z_0)\right\|^{1/\phi^n},
\]
for any norm $\|\cdot\| : \C^2\to \R$. Where defined, this $h_{f,z_0}$ satisfies the functional equation
\[
h_{f,z_0}(\mc{S}_f(x,y)) = \left[h_{f,z_0}(x,y)\right]^\phi.
\]
This is similar to the Böttcher functional equation $\Phi(g(z)) = [\Phi(z)]^n$ in the one-dimensional case, where $g(z) = z^n(a_n + a_{n+1}z + \cdots)$. Böttcher's theorem (cf. \cite{BuffEpsteinKochBottcherCoords}) says that the functional equation has a solution $\Phi$ which defines a holomorphic germ at the superattracting fixed point of $g$. Although our $h_{f,z_0}$ cannot possibly define an analogous holomorphic germ near $(z_0,z_0)$ (due in part to the irrational power in the functional equation and the fact that $h_{f,z_0}(z_0,z_0) = 0$), we shall prove that there is a holomorphic germ $\mc{H} : \mc{A}(f,z_0) \dashrightarrow \C$ near $(z_0,z_0)$ such that
\[
h_{f,z_0}(x,y) = |\mc{H}(x,y)|^{1/\phi}\cdot |x-z_0|^{1/\sqrt5} \cdot |y-z_0|^{\phi^{-1}/\sqrt5}.
\]
The germ $\mc{H}$ will have many of the nice properties enjoyed by the one-dimensional Böttcher coordinate $\Phi$. Therefore, we call it the Böttcher-type holomorphic germ for the secant map near $(z_0,z_0)$.

We remark that unlike other sources in the literature on the secant method dynamical system (e.g., \cite{BedfordFriggeSecant} and \cite{GarijoJarqueRealSecant}), we are not pre-supposing that $f$ is a polynomial or even a rational map. The potential function $h_{f,z_0}$ will always induce a holomorphic germ $\mc{H}$ regardless of the type of map $f$ is.

%%%%
\section{Power series and iterates of $\mc{S}_f$}
\label{sec:SecantIterates}

Before we can determine the iterates of $\mc{S}_f$ in the basin $\mc{A}(f,z_0)$, we require its power series about $(z_0,z_0)$. 

\begin{center}
\begin{minipage}{37em}
\noindent\textbf{Lemma \hypertarget{2:Jacobian}{2.1} (Jacobian)} \textit{If $z_0\in U$ is a simple root of $f$ (i.e., $f(z_0) = 0$ but $f'(z_0) \neq 0$), then the complex Jacobian of $\mc{S}_f$ is given by}
\[
D\mc{S}_f(x,y) = \begin{cases}
    \begin{bmatrix}
    -\frac{f(y)\left(f(x) - f(y) - f'(x)(x-y)\right)}{[f(x)-f(y)]^2} & \frac{f(x)\left(f(x) - f(y) - f'(y)(x-y)\right)}{[f(x)-f(y)]^2}\\
    1 & 0
    \end{bmatrix} & \text{if } x\neq y,\\[1.0em]
    \begin{bmatrix}
    \frac{f(x)f''(x)}{2[f'(x)]^2} & \frac{f(x)f''(x)}{2[f'(x)]^2}\\
    1 & 0
    \end{bmatrix} & \text{if } x = y
\end{cases}
\]
\textit{for all $(x,y)$ near $(z_0,z_0)$.}
\end{minipage}
\end{center}
\begin{proof}
We need the following two facts: first, let $a, b : U\to \C$ be holomorphic functions. Then the map $T_{a,b} : U\times U\to \C$ given by
\[
T_{a,b}(x,y) = \begin{cases}
    \dfrac{a(x)b(y) - a(y)b(x)}{x-y} & \text{if } x\neq y\\
    a'(x)b(x) - a(x)b'(x) & \text{if } x = y
\end{cases}
\]
is holomorphic. In particular,
the partials of $T_{a,b}$ along the diagonal are given by
\[
\frac{\partial T_{a,b}}{\partial x}(z,z) = \frac12\left[a''(z)b(z) - a(z)b''(z)\right] = \frac{\partial T_{a,b}}{\partial y}(z,z)
\]
for all $z\in U$. Second, in the notation of the first fact, we have
\[
\mc{S}_f(x,y) = \left(\frac{T_{f,\id}(x,y)}{T_{f,1}(x,y)}, x\right)
\]
for all $(x,y)$ sufficiently near $(z_0,z_0)$.

Let $(x_0,y_0)$ near $(z_0,z_0)$ be given. By our facts, we find
\begin{align*}
\frac{\partial \mc{S}_f}{\partial x}(x_0,y_0) & = \left.\frac{\partial}{\partial x}\left(\frac{T_{f,\id}(x,y_0)}{T_{f,1}(x,y_0)}, x\right)\right|_{x=x_0} = \left(\left.\frac{\partial}{\partial x}\left(\frac{T_{f,\id}(x,y_0)}{T_{f,1}(x,y_0)}\right)\right|_{x=x_0}, \left.\frac{\partial x}{\partial x}\right|_{x=x_0}\right)\\
& = \left(\frac{\frac{\partial T_{f,\id}}{\partial x}(x_0,y_0)\cdot T_{f,1}(x_0,y_0) - T_{f,\id}(x_0,y_0)\cdot \frac{\partial T_{f,1}}{\partial x}(x_0,y_0)}{\left[T_{f,1}(x_0,y_0)\right]^2}, 1\right).
\end{align*}
Likewise, 
\begin{align*}
\frac{\partial \mc{S}_f}{\partial y}(x_0,y_0) & = \left.\frac{\partial}{\partial y}\left(\frac{T_{f,\id}(x_0,y)}{T_{f,1}(x_0,y)}, x_0\right)\right|_{y=y_0} = \left(\left.\frac{\partial}{\partial y}\left(\frac{T_{f,\id}(x_0,y)}{T_{f,1}(x_0,y)}\right)\right|_{y=y_0}, \left.\frac{\partial x_0}{\partial y}\right|_{y=y_0}\right)\\
& = \left(\frac{\frac{\partial T_{f,\id}}{\partial y}(x_0,y_0)\cdot T_{f,1}(x_0,y_0) - T_{f,\id}(x_0,y_0)\cdot \frac{\partial T_{f,1}}{\partial y}(x_0,y_0)}{\left[T_{f,1}(x_0,y_0)\right]^2}, 0\right).
\end{align*}
If $x_0\neq y_0$, then the first coordinates evaluate to 
\begin{align*}
& \frac{\frac{\partial T_{f,\id}}{\partial x}(x_0,y_0)\cdot T_{f,1}(x_0,y_0) - T_{f,\id}(x_0,y_0)\cdot \frac{\partial T_{f,1}}{\partial x}(x_0,y_0)}{\left[T_{f,1}(x_0,y_0)\right]^2} = \left.\frac{\partial}{\partial x}\left(\frac{f(x)y_0 - f(y_0)x}{f(x) - f(y_0)}\right)\right|_{x=x_0}\\
& \ \ \ \ \ \ \ \ \ \ \ \ \ \ \ \ \ \ \ \ \ \ \ \ \ = \frac{(f'(x_0)y_0-f(y_0))(f(x_0)- f(y_0)) - (f(x_0)y_0 - f(y_0)x_0)f'(x_0)}{[f(x_0)-f(y_0)]^2}\\
& \ \ \ \ \ \ \ \ \ \ \ \ \ \ \ \ \ \ \ \ \ \ \ \ \ = -\frac{f(y_0)\left(f(x_0) - f(y_0) - f'(x_0)(x_0-y_0)\right)}{[f(x_0)-f(y_0)]^2}
\end{align*}
and
\begin{align*}
& \frac{\frac{\partial T_{f,\id}}{\partial y}(x_0,y_0)\cdot T_{f,1}(x_0,y_0) - T_{f,\id}(x_0,y_0)\cdot \frac{\partial T_{f,1}}{\partial y}(x_0,y_0)}{\left[T_{f,1}(x_0,y_0)\right]^2} = \left.\frac{\partial}{\partial y}\left(\frac{f(x_0)y - f(y)x_0}{f(x_0) - f(y)}\right)\right|_{y=y_0}\\
& \ \ \ \ \ \ \ \ \ \ \ \ \ \ \ \ \ \ \ \ \ \ \ \ \ = \frac{(f(x_0)-f'(y_0)x_0)(f(x_0)- f(y_0)) + (f(x_0)y_0 - f(y_0)x_0)f'(y_0)}{[f(x_0)-f(y_0)]^2}\\
& \ \ \ \ \ \ \ \ \ \ \ \ \ \ \ \ \ \ \ \ \ \ \ \ \ = \frac{f(x_0)\left(f(x_0) - f(y_0) - f'(y_0)(x_0-y_0)\right)}{[f(x_0)-f(y_0)]^2}.
\end{align*}
If instead $x_0 = w_0 = y_0$, then 
\begin{align*}
& \frac{\frac{\partial T_{f,\id}}{\partial x}(w_0,w_0)\cdot T_{f,1}(w_0,w_0) - T_{f,\id}(w_0,w_0)\cdot \frac{\partial T_{f,1}}{\partial x}(w_0,w_0)}{\left[T_{f,1}(w_0,w_0)\right]^2} \\
& \ \ \ \ \ \ \ \ \ \ \ \ \ \ \ \ \ \ \ \ \ \ \ \ \ =\frac{\frac{\partial T_{f,\id}}{\partial y}(w_0,w_0)\cdot T_{f,1}(w_0,w_0) - T_{f,\id}(w_0,w_0)\cdot \frac{\partial T_{f,1}}{\partial y}(w_0,w_0)}{\left[T_{f,1}(w_0,w_0)\right]^2}\\
& \ \ \ \ \ \ \ \ \ \ \ \ \ \ \ \ \ \ \ \ \ \ \ \ \ = \frac{\left(\frac12f''(w_0)w_0\right)f'(w_0) - \left(f'(w_0)w_0 - f(w_0)\right)\left(\frac12f''(w_0)\right)}{\left[f'(w_0)\right]^2}\\
& \ \ \ \ \ \ \ \ \ \ \ \ \ \ \ \ \ \ \ \ \ \ \ \ \ = \frac{f(w_0)f''(w_0)}{2[f'(w_0)]^2}.
\end{align*}
Thus, the complex Jacobian $D\mc{S}_f(x_0,y_0)$ is 
\[
D\mc{S}_f(x_0,y_0) = \begin{cases}
    \begin{bmatrix}
    -\frac{f(y_0)\left(f(x_0) - f(y_0) - f'(x_0)(x_0-y_0)\right)}{[f(x_0)-f(y_0)]^2} & \frac{f(x_0)\left(f(x_0) - f(y_0) - f'(y_0)(x_0-y_0)\right)}{[f(x_0)-f(y_0)]^2}\\
    1 & 0
    \end{bmatrix} & \text{if } x_0\neq y_0,\\[1.0em]
    \begin{bmatrix}
    \frac{f(x_0)f''(x_0)}{2[f'(x_0)]^2} & \frac{f(x_0)f''(x_0)}{2[f'(x_0)]^2}\\
    1 & 0
    \end{bmatrix} & \text{if } x_0 = y_0
\end{cases},
\]
as desired.
\end{proof}

\begin{center}
\begin{minipage}{37em}
\noindent\textbf{Lemma \hypertarget{2:PowerSeries}{2.2} (Power Series about $(z_0,z_0)$)} \textit{For each simple root $z_0\in U$ of $f$, there exists a unique holomorphic map $\mc{G}_{f,z_0} : \mc{A}(f,z_0) \to \C$ such that 
\[
\mc{S}_f(x,y) = \left(z_0 + \mc{G}_{f,z_0}(x,y)(x-z_0)(y-z_0), x\right)
\]
for all $(x,y)\in \mc{A}(f,z_0)$, and
\[
\mc{G}_{f,z_0}(z_0,z_0) = \frac{f''(z_0)}{2f'(z_0)}
\]
at the corresponding fixed point.}
\end{minipage}
\end{center}
\begin{proof}
\boxed{\textbf{Goal:}} Let $P_f\subseteq U$ be the set of poles of $f$ (it may be empty). Since $f$ is nonconstant, $P_f$ is a closed discrete subset of $U$. Thus, the set $(U\setminus P_f \times U\setminus P_f)\cap \mc{A}(f,z_0)$ is an open neighborhood of $(z_0,z_0)$, and there exists an $R > 0$ such that 
\[
(z_0,z_0) \in D_R(z_0) \times D_R(z_0) \subseteq (U\setminus P_f \times U\setminus P_f)\cap \mc{A}(f,z_0).
\] 
The bi-disk $D_R(z_0) \times D_R(z_0)$ will serve as a natural domain for the power series of $\mc{S}_f$ and $\mc{G}_{f,z_0}$.

For brevity in this proposition, we will let $\mc{S}_f^j$ (for $j\in\{1,2\}$) denote the coordinate functions of $\mc{S}_f$. Then the first coordinate $\mc{S}_f^1$ is expressible as the power series
\[
\mc{S}_f^1(x,y) = \sum_{d=0}^\infty \sum_{k=0}^d \frac1{(d-k)!\cdot k!}\cdot \frac{\partial^d \mc{S}_f^1}{\partial x^{d-k}\partial y^k}(z_0,z_0) \cdot(x-z_0)^{d-k}(y-z_0)^k,
\]
which converges for all $(x,y)\in D_R(z_0)\times D_R(z_0)$ since $\mc{S}_f^1$ is holomorphic there. We claim that the function $\mc{G}_{f,z_0} : \mc{A}(f,z_0)\to \C$ given by
\[
\mc{G}_{f,z_0}(x,y) = \begin{cases}
    \dfrac{\mc{S}_f^1(x,y) - z_0}{(x-z_0)(y-z_0)} & \text{if } x\neq z_0 \text{ and } y\neq z_0\\[0.7em]
    \dfrac1{x-z_0} - \dfrac{f'(z_0)}{f(x)}& \text{if } x\neq z_0 \text{ and } y = z_0\\[0.7em]
    \dfrac1{y-z_0} - \dfrac{f'(z_0)}{f(y)}& \text{if } x=z_0 \text{ and } y \neq z_0\\[0.7em]
    \displaystyle\sum_{d=2}^\infty \sum_{k=1}^{d-1} \frac1{(d-k)! \cdot k!}\cdot \frac{\partial^d \mc{S}_f^1}{\partial x^{d-k}\partial y^k}(z_0,z_0)\cdot (x-z_0)^{d-k-1}(y-z_0)^{k-1} & \text{if } x,y\in D_R(z_0)
\end{cases}
\]
is well-defined, holomorphic, and satisfies
\[
\mc{S}_f^1(x,y) = z_0 + \mc{G}_{f,z_0}(x,y)(x-z_0)(y-z_0)
\]
for all $(x,y) \in \mc{A}(f,z_0)$. 
\\

\noindent\boxed{\textbf{Step 1: Simplifying the Series:}} By Lemma \hyperlink{2:Jacobian}{2.1}, the first partials of $\mc{S}_f^1$ are given by
\[
\frac{\partial \mc{S}_f^1}{\partial x}(x,y) = \begin{cases}
    -\dfrac{f(y)(f(x)-f(y)-f'(x)(x-y))}{[f(x)-f(y)]^2} & \text{if } x\neq y\\
    \dfrac{f(x)f''(x)}{2[f'(x)]^2} & \text{if } x = y
\end{cases}
\]
and 
\[
\frac{\partial \mc{S}_f^1}{\partial y}(x,y) = \begin{cases}
    \dfrac{f(x)(f(x)-f(y)-f'(y)(x-y))}{[f(x)-f(y)]^2} & \text{if } x\neq y\\
    \dfrac{f(x)f''(x)}{2[f'(x)]^2} & \text{if } x = y
\end{cases}
\]
for all $(x,y)\in D_R(z_0) \times D_R(z_0)$. Because we need to evaluate at $(x,y) = (z_0,z_0)$, it will be useful to note that 
\[
\frac{\partial \mc{S}_f^1}{\partial x}(z_0,z_0) = 0 = \frac{\partial \mc{S}_f^1}{\partial y}(z_0,z_0)
\]
and for all $x,y\in D_R(z_0)\setminus\{z_0\}$, 
\[
\frac{\partial \mc{S}_f^1}{\partial x}(x,z_0) = 0, \ \ \frac{\partial \mc{S}_f^1}{\partial x}(z_0,y) = \frac{f(y) - f'(z_0)(y-z_0)}{f(y)}
\]
and
\[
\frac{\partial \mc{S}_f^1}{\partial y}(x,z_0) = \frac{f(x)-f'(z_0)(x-z_0)}{f(x)}, \ \ \frac{\partial \mc{S}_f^1}{\partial y}(z_0,y) = 0
\]
since $f(z_0) = 0$. 

We show via induction that the un-mixed partial derivatives all satisfy
\[
\frac{\partial^d \mc{S}_f^1}{\partial x^d}(x,z_0) = 0 = \frac{\partial^d \mc{S}_f^1}{\partial y^d}(z_0,y)
\]
for all $d\in\Z^+$ and all $x,y\in D_R(z_0)$. The base case $d = 1$ was noted above. Suppose the claim holds when $d = n$. Then
\[
\frac{\partial^{n+1} \mc{S}_f^1}{\partial x^{n+1}}(z_0,z_0) = \frac{d}{dx}\left.\left(\frac{\partial^n \mc{S}_f^1}{\partial x^n}(x,z_0)\right)\right|_{x=z_0} = \frac{d}{dx}\left.\left(0\right)\right|_{x=z_0} = 0
\]
and
\[
\frac{\partial^{n+1} \mc{S}_f^1}{\partial y^{n+1}}(z_0,z_0) = \frac{d}{dy}\left.\left(\frac{\partial^n \mc{S}_f^1}{\partial y^n}(z_0,y)\right)\right|_{x=z_0} = \frac{d}{dy}\left.\left(0\right)\right|_{y=z_0} = 0,
\]
as desired. 

Hence, by induction, we find
\[
\frac{\partial^d \mc{S}_f^1}{\partial x^d}(z_0,z_0) = 0 = \frac{\partial^d \mc{S}_f^1}{\partial y^d}(z_0,z_0)
\]
for all $d\in \Z^+$. The power series of $\mc{S}_f^1$ now simplifies to
\begin{align*}
\mc{S}_f^1(x,y) & = \sum_{d=0}^\infty \sum_{k=0}^d \frac1{(d-k)!\cdot k!}\cdot \frac{\partial^d \mc{S}_f^1}{\partial x^{d-k}\partial y^k}(z_0,z_0) \cdot(x-z_0)^{d-k}(y-z_0)^k\\
& = \frac1{0!\cdot 0!}\cdot \mc{S}_f^1(z_0,z_0)\cdot (x-z_0)^0(y-z_0)^0\\
& \ \ \ \ \ + \frac1{(1-0)!\cdot 0!}\cdot \frac{\partial \mc{S}_f^1}{\partial x}(z_0,z_0)(x-z_0) + \frac1{(1-1)!\cdot 1!}\cdot \frac{\partial \mc{S}_f^1}{\partial y}(z_0,z_0)(y-z_0)\\
& \ \ \ \ \ + \sum_{d=2}^\infty \sum_{k=0}^d \frac1{(d-k)!\cdot k!}\cdot \frac{\partial^d \mc{S}_f^1}{\partial x^{d-k}\partial y^k}(z_0,z_0) \cdot(x-z_0)^{d-k}(y-z_0)^k\\
& = z_0 + \sum_{d=2}^\infty \left[\frac1{(d-0)!\cdot 0!}\cdot \frac{\partial^d \mc{S}_f^1}{\partial x^d}(z_0,z_0) \cdot(x-z_0)^d + \frac1{(d-d)!\cdot d!}\cdot \frac{\partial^d \mc{S}_f^1}{\partial y^d}(z_0,z_0) \cdot(y-z_0)^d \right.\\
& \ \ \ \ \ \ \ \ \ \ \ \ \ \ \ \ \ \ \ \ \ \ \ \ \ \left. + \sum_{k=1}^{d-1} \frac1{(d-k)!\cdot k!}\cdot \frac{\partial^d \mc{S}_f^1}{\partial x^{d-k}\partial y^k}(z_0,z_0) \cdot(x-z_0)^{d-k}(y-z_0)^k\right]\\
& = z_0 + \sum_{d=2}^\infty \sum_{k=1}^{d-1} \frac1{(d-k)!\cdot k!}\cdot \frac{\partial^d \mc{S}_f^1}{\partial x^{d-k}\partial y^k}(z_0,z_0) \cdot(x-z_0)^{d-k}(y-z_0)^k\\
& = z_0 + (x-z_0)(y-z_0)\cdot \sum_{d=2}^\infty \sum_{k=1}^{d-1} \frac1{(d-k)!\cdot k!}\cdot \frac{\partial^d \mc{S}_f^1}{\partial x^{d-k}\partial y^k}(z_0,z_0) \cdot(x-z_0)^{d-1-k}(y-z_0)^{k-1}
\end{align*}
for all $x,y\in D_R(z_0)$. In particular, the final series on the far RHS must converge in the polydisk $D_R(z_0)\times D_R(z_0)$ since the original series does.
\\

\noindent\boxed{\textbf{Step 2: Well-Definedness:}} Observe that manipulating the equation at the end of Step 1 yields
\[
\frac{\mc{S}_f^1(x,y)-z_0}{(x-z_0)(y-z_0)} = \sum_{d=2}^\infty \sum_{k=1}^{d-1} \frac1{(d-k)!\cdot k!}\cdot \frac{\partial^d \mc{S}_f^1}{\partial x^{d-k}\partial y^k}(z_0,z_0) \cdot(x-z_0)^{d-1-k}(y-z_0)^{k-1}
\]
for all $(x,y) \in (D_R(z_0)\setminus\{z_0\}) \times (D_R(z_0)\setminus\{z_0\})$. So the first formula for $\mc{G}_{f,z_0}$ agrees with the fourth formula. We now show that the second and third formulas agree with the fourth.

When $y = z_0$, only the terms of the $\mc{G}_{f,z_0}$ series with $k = 1$ contribute:
\begin{align*}
&\left.\left(\sum_{d=2}^\infty \sum_{k=1}^{d-1} \frac1{(d-k)! \cdot k!}\cdot \frac{\partial^d \mc{S}_f^1}{\partial x^{d-k}\partial y^k}(z_0,z_0)\cdot (x-z_0)^{d-k-1}(y-z_0)^{k-1}\right)\right|_{y=z_0} \\
& \ \ \ \ \ = \sum_{d=2}^\infty  \frac1{(d-1)!}\cdot \frac{\partial^d \mc{S}_f^1}{\partial x^{d-1}\partial y}(z_0,z_0)\cdot (x-z_0)^{d-2} = \sum_{d=1}^\infty  \frac1{d!}\cdot \left.\frac{\partial^d}{\partial z^d}\left(\frac{\partial \mc{S}_f^1}{\partial y}(z,z_0)\right)\right|_{z=z_0}\cdot (x-z_0)^{d-1}\\
& \ \ \ \ \ = \frac1{x-z_0}\left[\frac{\partial \mc{S}_f^1}{\partial y}(x,z_0) - \frac{\partial \mc{S}_f^1}{\partial y}(z_0,z_0)\right] = \frac{f(x)-f'(z_0)(x-z_0)}{f(x)(x-z_0)} = \frac1{x-z_0} - \frac{f'(z_0)}{f(x)}
\end{align*}
as long as $x\in D_R(z_0)\setminus \{z_0\}$ (here we used the calculation in Step 1). Similarly, when $x = z_0$, only the terms with $k = d-1$ contribute:
\begin{align*}
&\left.\left(\sum_{d=2}^\infty \sum_{k=1}^{d-1} \frac1{(d-k)! \cdot k!}\cdot \frac{\partial^d \mc{S}_f^1}{\partial x^{d-k}\partial y^k}(z_0,z_0)\cdot (x-z_0)^{d-k-1}(y-z_0)^{k-1}\right)\right|_{x=z_0} \\
& \ \ \ \ \  = \sum_{d=2}^\infty  \frac1{(d-1)!}\cdot \frac{\partial^d \mc{S}_f^1}{\partial x\partial y^{d-1}}(z_0,z_0)\cdot (y-z_0)^{d-2} = \sum_{d=1}^\infty  \frac1{d!}\cdot \left.\frac{\partial^d}{\partial w^d}\left(\frac{\partial \mc{S}_f^1}{\partial x}(z_0,w)\right)\right|_{w=z_0}\cdot (y-z_0)^{d-1}\\
& \ \ \ \ \  = \frac1{y-z_0}\left[\frac{\partial \mc{S}_f^1}{\partial x}(z_0,y)\right] = \frac{f(y)-f'(z_0)(y-z_0)}{f(y)(y-z_0)} = \frac1{y-z_0} - \frac{f'(z_0)}{f(y)}
\end{align*}
as long as $y \in D_R(z_0)\setminus \{z_0\}$ (we used holomorphicity of $\mc{S}_f^1$ to switch the order of the partials).

Having checked all the possible domain overlap regions, we conclude the given formula for $\mc{G}_{f,z_0}$ is well-defined. So we have described a map $\mc{A}(f,z_0) \to \C$.  
\\

\noindent\boxed{\textbf{Step 3: Holomorphicity:}} Observe that $\mc{G}_{f,z_0}$ is holomorphic in the power series' polydisk $D_R(z_0)\times D_R(z_0)$ since the series converges there. Holomorphicity of $\mc{S}_f^1$ and closedness of the lines $\C\times \{z_0\}$ and $\{z_0\}\times \C$ in $\C^2$ also ensures that $\mc{G}_{f,z_0}$ is holomorphic at every point $(x,y)\in \mc{A}(f,z_0)$ such that $x\neq z_0$ and $y\neq z_0$. So it remains to check holomorphicity at points in the domain of the form $(x_0,z_0)$ and $(z_0,y_0)$, where $x_0,y_0\in \C\setminus D_R(z_0)$.

We shall use Hartogs' theorem on separate holomorphicity. Note that we have
\[
\mc{G}_{f,z_0}(x,z_0) = \frac1{x-z_0} - \frac{f'(z_0)}{f(x)} \text{ and } \mc{G}_{f,z_0}(z_0,y) = \frac1{y-z_0} - \frac{f'(z_0)}{f(y)}
\]
whenever $x,y\in \C\setminus \{z_0\}$. Therefore, since $f$ is meromorphic, we find $\mc{G}_{f,z_0}$ is separately holomorphic in the first argument at $(x_0,z_0)$ and in the second argument at $(z_0,y_0)$ (note that if $x_0$ or $y_0$ is a pole, then $z\mapsto 1/f(z)$ extends to be holomorphic in some open neighborhood of that pole). The corresponding partial derivatives at $(x_0,z_0)$ and $(z_0,y_0)$ are given by
\[
\frac{\partial \mc{G}_{f,z_0}}{\partial x}(x_0,z_0) = \left.\frac{d}{dx}\left(\mc{G}_{f,z_0}(x,z_0)\right)\right|_{x=x_0} = \left.\frac{d}{dx}\left(\frac1{x-z_0} - \frac{f'(z_0)}{f(x)}\right)\right|_{x=x_0} = \frac1{(x_0-z_0)^2} + \frac{f'(x_0)f'(z_0)}{[f(x_0)]^2}
\]
and similarly
\[
\frac{\partial \mc{G}_{f,z_0}}{\partial y}(z_0,y_0) = \frac1{(y_0-z_0)^2} + \frac{f'(y_0)f'(z_0)}{[f(y_0)]^2}.
\]

On the other hand, we have the formulas
\[
\mc{G}_{f,z_0}(x_0,y) = \begin{cases}
    \dfrac{\mc{S}_f^1(x_0,y) - z_0}{(x_0-z_0)(y-z_0)} & \text{if } y\neq z_0\\[0.7em]
    \dfrac1{x_0-z_0} - \dfrac{f'(z_0)}{f(x_0)} & \text{if } y = z_0
\end{cases} \text{ and } \mc{G}_{f,z_0}(x,y_0) = \begin{cases}
    \dfrac{\mc{S}_f^1(x,y_0) - z_0}{(x-z_0)(y_0-z_0)} & \text{if } x\neq z_0\\[0.7em]
    \dfrac1{y_0-z_0} - \dfrac{f'(z_0)}{f(y_0)} & \text{if } x = z_0
\end{cases}.
\]
The former is holomorphic in $y$ at least away from $y = y_0$, and the latter is holomorphic in $x$ away from $x = x_0$, because $\mc{S}_f^1$ is holomorphic. So it remains to check differentiability in $y$ (resp. $x$) at $y = z_0$ (resp. $x = z_0$). There are a few cases based on whether $x_0$ and/or $y_0$ is a pole of $f$.

If $x_0$ is a pole of $f$, then the first formula simplifies to
\begin{align*}
\mc{G}_{f,z_0}(x_0,y) & = \begin{cases}
    \dfrac{\mc{S}_f^1(x_0,y) - z_0}{(x_0-z_0)(y-z_0)} & \text{if } y\neq z_0\\[0.7em]
    \dfrac1{x_0-z_0} - \dfrac{f'(z_0)}{f(x_0)} & \text{if } y = z_0
\end{cases} = \begin{cases}
    \dfrac{y - z_0}{(x_0-z_0)(y-z_0)} & \text{if } y\in D_R(z_0)\setminus\{z_0\}\\[0.7em]
    \dfrac1{x_0-z_0} & \text{if } y = z_0
\end{cases} \\
& = \frac1{x_0-z_0},
\end{align*}
whence $\frac{\partial \mc{G}_{f,z_0}}{\partial y}(x_0,z_0) = 0$. Similarly, if $y_0$ is a pole of $f$, then the second formula becomes
\begin{align*}
\mc{G}_{f,z_0}(x,y_0) & = \begin{cases}
    \dfrac{\mc{S}_f^1(x,y_0) - z_0}{(x-z_0)(y_0-z_0)} & \text{if } x\neq z_0\\[0.7em]
    \dfrac1{y_0-z_0} - \dfrac{f'(z_0)}{f(y_0)} & \text{if } x = z_0
\end{cases} = \begin{cases}
    \dfrac{x - z_0}{(x-z_0)(y_0-z_0)} & \text{if } x\in D_R(z_0)\setminus\{z_0\}\\[0.7em]
    \dfrac1{y_0-z_0} & \text{if } x = z_0
\end{cases} \\
& = \frac1{y_0-z_0},
\end{align*}
whence $\frac{\partial \mc{G}_{f,z_0}}{\partial x}(z_0,y_0) = 0$. Now suppose $x_0$ and $y_0$ are both non-poles of $f$. Then the first formula is
\[
\mc{G}_{f,z_0}(x_0,y) = \begin{cases}
    \dfrac{\mc{S}_f^1(x_0,y) - z_0}{(x_0-z_0)(y-z_0)} & \text{if } y\neq z_0\\[0.7em]
    \dfrac1{x_0-z_0} - \dfrac{f'(z_0)}{f(x_0)} & \text{if } y = z_0
\end{cases} = \begin{cases}
    \dfrac{\frac{f(x_0)y-f(y)x_0}{f(x_0)-f(y)} - z_0}{(x_0-z_0)(y-z_0)} & \text{if } y\in D_R(z_0)\setminus\{z_0\}\\[0.7em]
    \dfrac{f(x_0)-f'(z_0)(x_0-z_0)}{f(x_0)(x_0-z_0)} & \text{if } y = z_0
\end{cases}
\]
and the second is 
\[
\mc{G}_{f,z_0}(x,y_0) = \begin{cases}
    \dfrac{\mc{S}_f^1(x,y_0) - z_0}{(x-z_0)(y_0-z_0)} & \text{if } x\neq z_0\\[0.7em]
    \dfrac1{y_0-z_0} - \dfrac{f'(z_0)}{f(y_0)} & \text{if } x = z_0
\end{cases} = \begin{cases}
    \dfrac{\frac{f(x)y_0-f(y_0)x}{f(x)-f(y_0)} - z_0}{(x-z_0)(y_0-z_0)} & \text{if } x\in D_R(z_0)\setminus\{z_0\}\\[0.7em]
    \dfrac{f(y_0)-f'(z_0)(y_0-z_0)}{f(y_0)(y_0-z_0)} & \text{if } x = z_0
\end{cases}.
\]
At this point, it will be useful to use the power series of $f(z)$ in $D_R(z_0)$ about $z = z_0$ to say there is a (unique) holomorphic map $g_{z_0} : D_R(z_0) \to \C$ satisfying
\[
f(z) = f'(z_0)(z-z_0) + g_{z_0}(z)(z-z_0)^2 \text{ and } g_{z_0}(z_0) = \frac12f''(z_0)
\]
for all $z\in D_R(z_0)$. 

Observe now that for any $w_0\in U\setminus P_f$, and all $z\in D_R(z_0)\setminus\{z_0\}$ with $f(z) \neq f(w_0)$, we have 
\begin{align*}
\frac{f(z)w_0 - f(w_0)z}{f(z) - f(w_0)} - z_0 
& = \frac{\left(f(z)w_0 - f(w_0)z\right) - z_0\left(f(z)-f(w_0)\right)}{f(z)-f(w_0)} = \frac{f(z)(w_0-z_0) - f(w_0)(z-z_0)}{f(z)-f(w_0)}\\
& = \frac{\left(f'(z_0)(w_0-z_0) - f(w_0)\right)(z-z_0) + g_{z_0}(z)(w_0-z_0)(z-z_0)^2}{f(z)-f(w_0)}\\
& = \frac{f'(z_0)(w_0-z_0) - f(w_0) + g_{z_0}(z)(w_0-z_0)(z-z_0)}{f(z)-f(w_0)}\cdot (z-z_0).
\end{align*}
Therefore,
\begin{align*}
& \frac{\partial\mc{G}_{f,z_0}}{\partial y}(x_0,z_0) \\
& \ \ \ \ \ \ \ = \lim_{y\to z_0} \frac{\mc{G}_{f,z_0}(x_0,y) - \mc{G}_{f,z_0}(x_0,z_0)}{y-z_0} \\
& \ \ \ \ \ \ \ = \lim_{y\to z_0} \frac1{y-z_0}\left[\frac{\frac{f(x_0)y-f(y)x_0}{f(x_0)-f(y)} - z_0}{(x_0-z_0)(y-z_0)} - \frac{f(x_0)-f'(z_0)(x_0-z_0)}{f(x_0)(x_0-z_0)}\right]\\
& \ \ \ \ \ \ \ = \lim_{y\to z_0} \frac1{y-z_0}\left[\frac{f'(z_0)(x_0-z_0) - f(x_0) + g_{z_0}(y)(x_0-z_0)(y-z_0)}{(f(y)-f(x_0))(x_0-z_0)} - \frac{f(x_0)-f'(z_0)(x_0-z_0)}{f(x_0)(x_0-z_0)}\right]\\
& \ \ \ \ \ \ \ = \left.\frac{d}{dy}\left(\frac{f'(z_0)(x_0-z_0) - f(x_0) + g_{z_0}(y)(x_0-z_0)(y-z_0)}{(f(y)-f(x_0))(x_0-z_0)}\right)\right|_{y=z_0} \\
& \ \ \ \ \ \ \ = \frac{g_{z_0}(z_0)(x_0-z_0)(f(z_0)-f(x_0)) - \left(f'(z_0)(x_0-z_0) - f(x_0)\right)f'(z_0)}{\left[f(z_0)-f(x_0)\right]^2(x_0-z_0)}\\
& \ \ \ \ \ \ \ = \frac{-\frac12f''(z_0)f(x_0)(x_0-z_0) - \left(f'(z_0)(x_0-z_0) - f(x_0)\right)f'(z_0)}{[f(x_0)]^2\cdot (x_0-z_0)}.
\end{align*}
Similarly,
\begin{align*}
& \frac{\partial\mc{G}_{f,z_0}}{\partial x}(z_0,y_0) \\
& \ \ \ \ \ \ \ = \lim_{x\to z_0} \frac{\mc{G}_{f,z_0}(x,y_0) - \mc{G}_{f,z_0}(z_0,y_0)}{x-z_0} \\
& \ \ \ \ \ \ \ = \lim_{x\to z_0} \frac1{x-z_0}\left[\frac{\frac{f(x)y_0-f(y_0)x}{f(x)-f(y_0)} - z_0}{(x-z_0)(y_0-z_0)} - \frac{f(x)-f'(z_0)(x-z_0)}{f(x)(x-z_0)}\right]\\
& \ \ \ \ \ \ \ = \lim_{x\to z_0} \frac1{x-z_0}\left[\frac{f'(z_0)(y_0-z_0) - f(y_0) + g_{z_0}(x)(y_0-z_0)(x-z_0)}{(f(x)-f(y_0))(y_0-z_0)} - \frac{f(y_0)-f'(z_0)(y_0-z_0)}{f(y_0)(y_0-z_0)}\right]\\
& \ \ \ \ \ \ \ = \left.\frac{d}{dx}\left(\frac{f'(z_0)(y_0-z_0) - f(y_0) + g_{z_0}(x)(y_0-z_0)(x-z_0)}{(f(x)-f(y_0))(y_0-z_0)}\right)\right|_{x=z_0} \\
& \ \ \ \ \ \ \ = \frac{g_{z_0}(z_0)(y_0-z_0)(f(z_0)-f(y_0)) - \left(f'(z_0)(y_0-z_0) - f(y_0)\right)f'(z_0)}{\left[f(z_0)-f(y_0)\right]^2(y_0-z_0)}\\
& \ \ \ \ \ \ \ = \frac{-\frac12f''(z_0)f(y_0)(y_0-z_0) - \left(f'(z_0)(y_0-z_0) - f(y_0)\right)f'(z_0)}{[f(y_0)]^2\cdot (y_0-z_0)}.
\end{align*}
We conclude in all cases that $\mc{G}_{f,z_0}$ is separately holomorphic in the second argument at $(x_0,z_0)$ and in the first argument at $(z_0,y_0)$. 

So $\mc{G}_{f,z_0}$ is separately holomorphic at all points of the form $(x_0,z_0)$ and $(z_0,y_0)$ for certain $x_0,y_0\in \C\setminus D_R(z_0)$. Since $\mc{G}_{f,z_0}$ is already known to be fully holomorphic away from these points, it follows from Hartogs' theorem that $\mc{G}_{f,z_0}$ is fully holomorphic on all of $\mc{A}(f,z_0)$. 
\\

\noindent\boxed{\textbf{Step 4: Other Properties of $\mc{G}_{f,z_0}$:}} By definition of $\mc{G}_{f,z_0}$, we have
\[
\mc{S}_f^1(x,y) = z_0 + \mc{G}_{f,z_0}(x,y)(x-z_0)(y-z_0)
\]
for all $(x,y)\in \mc{A}(f,z_0)$ such that $x\neq z_0$ and $y\neq z_0$. In particular, this implies
\[
\mc{S}_f(x,y) = \left(z_0 + \mc{G}_{f,z_0}(x,y)(x-z_0)(y-z_0), x\right)
\]
as required. Since the domain $\mc{A}(f,z_0)$ is open in $\C^2$ and the lines $\{z_0\}\times \C$ and $\C\times \{z_0\}$ are both closed in $\C^2$, it follows from the Uniqueness Principle for holomorphic maps that $\mc{G}_{f,z_0}$ must be unique in this respect.

It remains to compute $\mc{G}_{f,z_0}(z_0,z_0)$. Referring back to the power series definition of $\mc{G}_{f,z_0}$ near $(z_0,z_0)$, we find that only the $d = 2$, $k= 1$ term contributes:
\[
\mc{G}_{f,z_0}(z_0,z_0) = \frac1{(2-1)!\cdot 1!}\cdot \frac{\partial^2 \mc{S}_f^1}{\partial x^{2-1}\partial y^1}(z_0,z_0) = \frac{\partial^2 \mc{S}_f^1}{\partial x\partial y}(z_0,z_0).
\]
Now the calculations in Steps 1 and 3 can be used to obtain
\begin{align*}
\frac{\partial^2 \mc{S}_f^1}{\partial x\partial y}(z_0,z_0) &= \left.\frac{\partial}{\partial x}\left(\frac{\partial \mc{S}_f^1}{\partial y}(x,z_0)\right)\right|_{x=z_0} = \lim_{x\to z_0} \frac{\frac{\partial \mc{S}_f^1}{\partial y}(x,z_0) - \frac{\partial \mc{S}_f^1}{\partial y}(z_0,z_0)}{x-z_0}\\
& = \lim_{x\to z_0} \frac1{x-z_0}\cdot\left[\frac{f(x)-f'(z_0)(x-z_0)}{f(x)} - 0\right] = \lim_{x\to z_0} \frac{g_{z_0}(x)(x-z_0)^2}{f'(z_0)(x-z_0)^2 + g_{z_0}(z)(x-z_0)^3}\\
& = \lim_{x\to z_0} \frac{g_{z_0}(x)}{f'(z_0) + g_{z_0}(x)(x-z_0)} = \frac{\frac12f''(z_0)}{f'(z_0)}\\
& = \frac{f''(z_0)}{2f'(z_0)},
\end{align*}
as desired.
\end{proof}

The remarkable property of the formula for $\mc{S}_f$ in Lemma \hyperlink{2:PowerSeries}{2.2} is that $(x-z_0)(y-z_0)$ can be factored out of all of the degree-$2$ and higher terms in the power series. It is this fact that allows us to ensure the potential $h_{f,z_0}$ induces a holomorphic germ $\mc{H}$, not simply a real-analytic modulus.

Observe that if $f''(z_0) = 0$, then $\mc{G}_{f,z_0}(z_0,z_0) = 0$, and the term $\mc{G}_{f,z_0}(x,y)(x-z_0)(y-z_0)$ in the power series of $\mc{S}_f$ does not actually begin with a quadratic term, rather, it will begin with a higher degree term. In such cases, one finds that the convergence of $\mc{S}_f^{\circ n}(x,y)$ will actually be faster than the aforementioned rate $\phi$, that is, the potential $h_{f,z_0}$ will be identically zero. We shall hereafter avoid this exceptional case, and assume $f''(z_0) \neq 0$. 

\begin{center}
\begin{minipage}{37em}
\noindent\textbf{Proposition \hypertarget{2:Iterates}{2.3} (Iterates on Basin of Attraction)} \textit{On the basin of attraction $\mc{A}(f,z_0)$ of $(z_0,z_0)$, the iterates $\mc{S}_f^{\circ n}$ take the form
    \begin{align*}
    \mc{S}_f^{\circ n}(x,y) & = \left(z_0 + \mc{G}_n(x,y)(x-z_0)^{F_{n+1}}(y-z_0)^{F_n},\right.\\ 
    & \ \ \ \ \  \ \ \left. z_0 + \mc{G}_{n-1}(x,y)(x-z_0)^{F_n}(y-z_0)^{F_{n-1}}\right)
    \end{align*}
    for all $n\in\Z^+$ and $(x,y)\in \mc{A}(f,z_0)$. Here $F_n$ denotes the $n$th Fibonacci number, determined by $F_{n+1} = F_n + F_{n-1}$ and $F_0 = 0$ and $F_1 = 1$. The functions $\mc{G}_n : \mc{A}(f,z_0) \to \C$ are given by
    \[
    \mc{G}_0(x,y) = 1 \text{ and } \mc{G}_n(x,y) = \prod_{k=0}^{n-1} \left[\mc{G}_{f,z_0}(\mc{S}_f^{\circ k}(x,y))\right]^{F_{n-k}}
    \]
    for $n\in\Z^+$, where $\mc{G}_{f,z_0} : \mc{A}(f,z_0) \to \C$ is the function from Lemma \hyperlink{2:PowerSeries}{2.2}.}
\end{minipage}
\end{center}
\begin{proof}
We prove this claim by induction on $n\in \Z^+$. The base case $n = 1$ is the result of Lemma \hyperlink{2:PowerSeries}{2.2}: we have
\begin{align*}
\mc{S}_f(x,y) & = \left(z_0 + \mc{G}_{f,z_0}(x,y)(x-z_0)(y-z_0), x\right)\\
& = \left(z_0 + \left(\prod_{k=0}^{1-1} \left[\mc{G}_{f,z_0}(\mc{S}_f^{\circ k}(x,y))\right]^{F_{1-0}}\right)\cdot (x-z_0)^1(y-z_0)^1, z_0 + 1\cdot (x-z_0)^1(y-z_0)^0\right)\\
& = \left(z_0 + \mc{G}_1(x,y)(x-z_0)^{F_{1+1}}(y-z_0)^{F_1}, z_0 + \mc{G}_0(x,y)(x-z_0)^{F_1}(y-z_0)^{F_0}\right)
\end{align*}
for all $(x,y)\in \mc{A}(f,z_0)$. For the inductive step, suppose the claim holds when $n = m$, for some $m\in \Z^+$. Then
\begin{align*}
\mc{S}_f^{\circ (m+1)}(x,y) & = \mc{S}_f^{\circ m}\left(\mc{S}_f(x,y)\right) = \mc{S}_f^{\circ m}\left(z_0 + \mc{G}_{f,z_0}(x,y)(x-z_0)(y-z_0), x\right)\\
& = \left(z_0 + \mc{G}_m(\mc{S}_f(x,y))\cdot \left[\mc{G}_{f,z_0}(x,y)(x-z_0)(y-z_0)\right]^{F_{m+1}}\cdot (x-z_0)^{F_m},\right.\\
& \ \ \ \ \ \ \ \left. z_0 + \mc{G}_{m-1}(\mc{S}_f(x,y))\cdot \left[\mc{G}_{f,z_0}(x,y)(x-z_0)(y-z_0)\right]^{F_m}\cdot (x-z_0)^{F_{m-1}}\right)\\
& = \left(z_0 + \left(\prod_{k=0}^{m-1} \left[\mc{G}_{f,z_0}\left(\mc{S}_f^{\circ k}(\mc{S}_f(x,y))\right)\right]^{F_{m-k}}\right)\cdot [\mc{G}_{f,z_0}(x,y)]^{F_{m+1}} \cdot (x-z_0)^{F_{m+1}+F_m}(y-z_0)^{F_{m+1}},\right.\\
& \ \ \ \ \ \ \ \ \  \left. z_0 + \left(\prod_{k=0}^{m-2} \left[\mc{G}_{f,z_0}\left(\mc{S}_f^{\circ k}(\mc{S}_f(x,y))\right)\right]^{F_{m-1-k}}\right)\cdot [\mc{G}_{f,z_0}(x,y)]^{F_m} \cdot (x-z_0)^{F_m+F_{m-1}}(y-z_0)^{F_m}\right)\\
& = \left(z_0 + \left(\prod_{k=1}^m \left[\mc{G}_{f,z_0}(x,y)\right]^{F_{m+1-k}}\right)\cdot [\mc{G}_{f,z_0}(x,y)]^{F_{m+1}} \cdot (x-z_0)^{F_{m+2}}(y-z_0)^{F_{m+1}},\right.\\
& \ \ \ \ \ \ \ \ \ \ \ \ \ \ \  \left. z_0 + \left(\prod_{k=1}^m \left[\mc{G}_{f,z_0}(x,y)\right]^{F_{m-k}}\right)\cdot [\mc{G}_{f,z_0}(x,y)]^{F_m} \cdot (x-z_0)^{F_{m+1}}(y-z_0)^{F_m}\right)\\
& = \left(z_0 + \left(\prod_{k=0}^m \left[\mc{G}_{f,z_0}(x,y)\right]^{F_{m+1-k}}\right)\cdot (x-z_0)^{F_{m+2}}(y-z_0)^{F_{m+1}}, \right.\\
& \ \ \ \ \ \ \ \ \ \ \ \  \left.z_0 +\left(\prod_{k=0}^{m-1} \left[\mc{G}_{f,z_0}(x,y)\right]^{F_{m-k}}\right)\cdot (x-z_0)^{F_{m+1}}(y-z_0)^{F_m}\right)\\
& = \left(z_0 + \mc{G}_{m+1}(x,y)(x-z_0)^{F_{m+2}}(y-z_0)^{F_{m+1}}, z_0 + \mc{G}_m(x,y)(x-z_0)^{F_{m+1}}(y-z_0)^{F_m}\right)
\end{align*}
for all $(x,y)\in \mc{A}(f,z_0)$, because then $\mc{S}_f(x,y)\in \mc{A}(f,z_0)$ too. By induction, the claim holds for all $n\in \Z^+$.
\end{proof}

%%%%
\section{The Böttcher-type construction}
Under the assumption that $f''(z_0) \neq 0$, we may take any $w_0$ in the countable set $\log\left(\frac{f''(z_0)}{2f'(z_0)}\right)$. Each logarithm will correspond in a canonical way to a different holomorphic germ $\mc{H}_{f,z_0,w_0} : \mc{A}(f,z_0) \dashrightarrow \C$ induced by the potential $h_{f,z_0}$. This is similar to how the one-dimensional Böttcher coordinate is unique up to multiplication by a root of unity. All of the germs $\mc{H}_{f,z_0,w_0}$ we will obtain will have the same modulus, and therefore we shall find this single modulus will extend to be continuous on all of $\mc{A}(f,z_0)$.
\newpage

\begin{center}
\begin{minipage}{37em}
\noindent\textbf{Theorem \hypertarget{3:BöttcherFunction}{3.1} (Böttcher-Type Functions)} \textit{Suppose we have $f''(z_0) \neq 0$.
\begin{itemize}
    \item[\textbf{(1).}] For each choice of logarithm $w_0\in \log\left(\frac{f''(z_0)}{2f'(z_0)}\right)$, there exists an open set $V\subseteq \C^2$ satisfying $(z_0,z_0)\in V\subseteq \mc{A}(f,z_0)$ and $\mc{S}_f(V)\subseteq V$, and there exists a holomorphic map $\mc{H}_{f,z_0,w_0} : V\to \C$ with
    \[
    \mc{H}_{f,z_0,w_0}(z_0,z_0) = \exp\left(\frac{5+3\sqrt{5}}{10}\cdot w_0\right)=: \left(\frac{f''(z_0)}{2f'(z_0)}\right)^{\frac{5+3\sqrt{5}}{10}} 
    \]
    and
    \[
    \lim_{n\to \infty} \left|(\pi_1\circ \mc{S}_f^{\circ n})(x,y) - z_0\right|^{1/\phi^n} = |\mc{H}_{f,z_0,w_0}(x,y)|\cdot |x-z_0|^{\phi/\sqrt{5}}\cdot |y-z_0|^{1/\sqrt5}
    \]
    for all $(x,y)\in V$, where $\pi_1 : \C^2\to \C$ is the first coordinate projection, and $\phi = \frac12\left(1+\sqrt{5}\right)$ is the golden ratio.
    \item[\textbf{(2).}] For each fixed $w_0 \in \log\left(\frac{f''(z_0)}{2f'(z_0)}\right)$, the corresponding $\mc{H}_{f,z_0,w_0}$ is unique in the sense that any two such functions with the same properties agree on the overlap of their domains. 
    \item[\textbf{(3).}] %Every possible function $\mc{H}_{f,z_0,w_0}$ has the same modulus where their domains overlap, regardless of initial logarithm chosen. There
    There exists a unique continuous function $\widehat{\mc{H}}_{f,z_0} : \mc{A}(f,z_0) \to \R$ such that
    \[
    \lim_{n\to \infty} \left|(\pi_1\circ \mc{S}_f^{\circ n})(x,y) - z_0\right|^{1/\phi^n} = \widehat{\mc{H}}_{f,z_0}(x,y)\cdot |x-z_0|^{\phi/\sqrt{5}}\cdot |y-z_0|^{1/\sqrt5}
    \]
    for all $(x,y)\in \mc{A}(f,z_0)$, and
    \[
    \widehat{\mc{H}}_{f,z_0}(x,y) = |\mc{H}_{f,z_0,w_0}(x,y)|
    \]
    for every function $\mc{H}_{f,z_0,w_0} : V\to \C$ from (1) (even letting $w_0$ vary), and every $(x,y)\in V$. 
    \item[\textbf{(4).}] The map $\widehat{\mc{H}}_{f,z_0}$ from (3) satisfies
    \[
    \widehat{\mc{H}}_{f,z_0}(x,y) = \left[\widehat{\mc{H}}_{f,z_0}(\mc{S}_f(x,y))\right]^{1/\phi} \cdot |\mc{G}_{f,z_0}(x,y)|^{1/\sqrt5}
    \]
    for all $(x,y)\in \mc{A}(f,z_0)$, where $\mc{G}_{f,z_0} : \mc{A}(f,z_0) \to \C$ is as Lemma \hyperlink{2:PowerSeries}{2.2}. Therefore, $\widehat{\mc{H}}_{f,x_0}$ is real analytic wherever it is nonzero (i.e., on the preimage $\widehat{\mc{H}}_{f,z_0}^{-1}(\R\setminus\{0\})$).
\end{itemize}}
\end{minipage}
\end{center}
\begin{proof}
\boxed{\textbf{(1).}} \underline{\textbf{Step 1:}} Let $w_0\in \log\left(\frac{f''(z_0)}{2f'(z_0)}\right)$ be given. Let $\mc{G}_{f,z_0} : \mc{A}(f,z_0) \to \C$ be the function from Lemma \hyperlink{2:PowerSeries}{2.2}. Then we know
\[
\mc{G}_{f,z_0}(z_0,z_0) = \frac{f''(z_0)}{2f'(z_0)} = \exp(w_0).
\]
We also have from Proposition \hyperlink{2:Iterates}{2.3} that 
\[
(\pi_1\circ \mc{S}_f^{\circ n})(x,y) = z_0 + \mc{G}_n(x,y)(x-z_0)^{F_{n+1}}(y-z_0)^{F_n}
\]
for all $n\in \Z^+$ and all $(x,y)\in \mc{A}(f,z_0)$, where $\mc{G}_n : \mc{A}(f,z_0) \to \C$ is given by
\[
\mc{G}_n(x,y) = \prod_{k=0}^{n-1} \left[\mc{G}_{f,z_0}(\mc{S}_f^{\circ k}(x,y))\right]^{F_{n-k}}.
\]
Our goal is to explicitly construct holomorphic branches of the maps $(x,y) \mapsto [\mc{G}_n(x,y)]^{1/\phi^n}$ in some open neighborhood of $(z_0,z_0)$, then show our branches converge uniformly on that neighborhood as $n\to\infty$.

Since clearly $\mc{G}_{f,z_0}(z_0,z_0) \neq 0$ (because we assumed $f''(z_0)\neq 0$) and $\mc{G}_{f,z_0}$ is continuous, it follows that there exists a $\rho > 0$ such that 
\[
D_\rho(z_0)\times D_\rho(z_0) \subseteq (U\setminus P_f \times U\setminus P_f) \cap \mc{A}(f,z_0)
\]
(where $P_f \subseteq U$ is the polar locus of $f$) and we have
\[
(x,y)\in D_\rho(z_0)\times D_\rho(z_0) \text{ implies } \left|\mc{G}_{f,z_0}(x,y) - \frac{f''(z_0)}{2f'(z_0)}\right| < \frac12\cdot\left|\frac{f''(z_0)}{2f'(z_0)}\right|.
\]
In particular, for these $(x,y)$, we have
\[
|\mc{G}_{f,z_0}(x,y)| \leq \left|\mc{G}_{f,z_0}(x,y) - \frac{f''(z_0)}{2f'(z_0)}\right| + \left|\frac{f''(z_0)}{2f'(z_0)}\right| < \frac32\cdot\left|\frac{f''(z_0)}{2f'(z_0)}\right|.
\]
Fix 
\[
r = \min\left\{\rho, \frac23\cdot\left|\frac{2f'(z_0)}{f''(z_0)}\right|\right\} > 0.
\]
Then for all $(x,y)\in D_r(z_0)\times D_r(z_0)$, we have
\begin{align*}
\left|(\pi_1\circ \mc{S}_f)(x,y)-z_0\right| & = \left|\mc{G}_{f,z_0}(x,y)(x-z_0)(y-z_0)\right| =  |\mc{G}_{f,z_0}(x,y)|\cdot|x-z_0|\cdot|y-z_0|\\
& < \frac32\cdot\left|\frac{f''(z_0)}{2f'(z_0)}\right|\cdot r \cdot r \leq r.
\end{align*}
If $\pi_2 : \C^2\to\C$ is the second coordinate projection, then we also have
\[
\left|(\pi_2\circ\mc{S}_f)(x,y)-z_0\right| = |x-z_0| < r.
\]
Therefore, we conclude
\[
\mc{S}_f(D_r(z_0)\times D_r(z_0)) \subseteq D_r(z_0)\times D_r(z_0). 
\]
This bi-disk shall be our desired $V$.
\\

\noindent\underline{\textbf{Step 2:}} Observe that the disk $D_{\frac12\cdot\left|\frac{f''(z_0)}{2f'(z_0)}\right|}\left(\frac{f''(z_0)}{2f'(z_0)}\right)$ does not contain $0$. Because the disk is simply connected, there exists a holomorphic branch $\mc{L} : D_{\frac12\cdot\left|\frac{f''(z_0)}{2f'(z_0)}\right|}\left(\frac{f''(z_0)}{2f'(z_0)}\right) \to \C$ of the complex logarithm such that
\[
\mc{L}\left(\frac{f''(z_0)}{2f'(z_0)}\right) = w_0.
\]
Now for each $n\in \Z^+$, define a map $H_n : D_r(z_0)\times D_r(z_0) \to \C$ by
\[
H_n(x,y) = \exp\left(\sum_{k=0}^{n-1} \frac{F_{n-k}}{\phi^n}\cdot \mc{L}\left(\mc{G}_{f,z_0}(\mc{S}_f^{\circ k}(x,y))\right)\right).
\]
Then $H_n$ is holomorphic, because $\mc{S}_f^{\circ k}(D_r(z_0)\times D_r(z_0))\subseteq D_r(z_0)\times D_r(z_0)$ by induction on $k$, and certainly $\mc{L}\circ \mc{G}_{f,z_0}$ is holomorphic on $D_r(z_0)\times D_r(z_0)$ by definition of $r$ and $\rho$. 

The functions $H_n$ relate to the functions $\mc{G}_n$ in the following manner:
\begin{align*}
H_n(x,y) & = \exp\left(\sum_{k=0}^{n-1} \frac{F_{n-k}}{\phi^n}\cdot \mc{L}\left(\mc{G}_{f,z_0}(\mc{S}_f^{\circ k}(x,y))\right)\right) = \prod_{k=0}^{n-1}\exp\left( \frac{F_{n-k}}{\phi^n}\cdot \mc{L}\left(\mc{G}_{f,z_0}(\mc{S}_f^{\circ k}(x,y))\right)\right)\\
& = \prod_{k=0}^{n-1}\exp\left( \frac1{\phi^n}\cdot \mc{L}\left(\left[\mc{G}_{f,z_0}(\mc{S}_f^{\circ k}(x,y))\right]^{F_{n-k}}\right)\right) = \prod_{k=0}^{n-1} \left[\left[\mc{G}_{f,z_0}(\mc{S}_f^{\circ k}(x,y))\right]^{F_{n-k}}\right]^{1/\phi^n}\\
& = \left[\prod_{k=0}^{n-1}\left[\mc{G}_{f,z_0}(\mc{S}_f^{\circ k}(x,y))\right]^{F_{n-k}}\right]^{1/\phi^n} = \left[\mc{G}_n(x,y)\right]^{1/\phi^n}.
\end{align*}
Thus, the $H_n$ are the desired holomorphic branches mentioned in Step 1. We now prove that the $H_n$ converge uniformly on $D_r(z_0)\times D_r(z_0)$. 

%It is sufficient to show the $g_n$ are uniformly Cauchy on the polydisk, because $\C$ is complete. 

Since $\exp$ is holomorphic and $\C$ is complete, it is sufficient to prove that the sequence of functions
\[
\sum_{k=0}^{n-1} \frac{F_{n-k}}{\phi^n}\cdot \mc{L}\left(\mc{G}_{f,z_0}(\mc{S}_f^{\circ k}(x,y))\right)
\]
is uniformly Cauchy for $(x,y)\in D_r(z_0)\times D_r(z_0)$. Observe that for all $m,n\in\Z^+$ with $m > n$, we have
\begin{align*}
& \left|\sum_{k=0}^{m-1} \frac{F_{m-k}}{\phi^m}\cdot \mc{L}\left(\mc{G}_{f,z_0}(\mc{S}_f^{\circ k}(x,y))\right) - \sum_{k=0}^{n-1} \frac{F_{n-k}}{\phi^n}\cdot \mc{L}\left(\mc{G}_{f,z_0}(\mc{S}_f^{\circ k}(x,y))\right)\right|\\
& \ \ \ \ \ = \left|\sum_{j=n}^{m-1} \left[\sum_{k=0}^j \frac{F_{j+1-k}}{\phi^{j+1}}\cdot \mc{L}\left(\mc{G}_{f,z_0}(\mc{S}_f^{\circ k}(x,y))\right) - \sum_{k=0}^{j-1} \frac{F_{j-k}}{\phi^j}\cdot \mc{L}\left(\mc{G}_{f,z_0}(\mc{S}_f^{\circ k}(x,y))\right)\right]\right|\\
& \ \ \ \ \ = \left|\sum_{j=n}^{m-1}\left[\frac{F_1}{\phi^{j+1}}\cdot \mc{L}\left(\mc{G}_{f,z_0}(\mc{S}_f^{\circ j}(x,y))\right) + \sum_{k=0}^{j-1} \left(\frac{F_{j+1-k}}{\phi^{j+1}} -\frac{F_{j-k}}{\phi^j}\right)\cdot \mc{L}\left(\mc{G}_{f,z_0}(\mc{S}_f^{\circ k}(x,y))\right)\right]\right|\\
& \ \ \ \ \ \leq \sum_{j=n}^{m-1}\left|\frac{F_1}{\phi^{j+1}}\cdot \mc{L}\left(\mc{G}_{f,z_0}(\mc{S}_f^{\circ j}(x,y))\right) + \sum_{k=0}^{j-1} \left(\frac{F_{j+1-k}}{\phi^{j+1}} -\frac{F_{j-k}}{\phi^j}\right)\cdot \mc{L}\left(\mc{G}_{f,z_0}(\mc{S}_f^{\circ k}(x,y))\right)\right|\\
& \ \ \ \ \ \leq \sum_{j=n}^{m-1} \left[\frac{F_1}{\phi^{j+1}}\cdot \left|\mc{L}\left(\mc{G}_{f,z_0}(\mc{S}_f^{\circ j}(x,y))\right)\right| + \sum_{k=0}^{j-1} \left|\frac{F_{j+1-k}}{\phi^{j+1}} -\frac{F_{j-k}}{\phi^j}\right|\cdot \left|\mc{L}\left(\mc{G}_{f,z_0}(\mc{S}_f^{\circ k}(x,y))\right)\right|\right]\\
& \ \ \ \ \ \leq \sum_{j=n}^{m-1} \left[\left(\frac{F_1}{\phi^{j+1}} + \sum_{k=0}^{j-1} \left|\frac{F_{j+1-k}}{\phi^{j+1}} -\frac{F_{j-k}}{\phi^j}\right|\right)\cdot \sup\left\{\left|\mc{L}\left(\mc{G}_{f,z_0}(a,b)\right)\right| : (a,b)\in D_r(z_0)\times D_r(z_0)\right\}\right]\\
& \ \ \ \ \ = \left(\sum_{j=n}^{m-1} \left[\frac{F_1}{\phi^{j+1}} + \sum_{k=0}^{j-1} \left|\frac{F_{j+1-k}}{\phi^{j+1}} -\frac{F_{j-k}}{\phi^j}\right|\right]\right)\cdot \sup\left\{\left|\mc{L}\left(\mc{G}_{f,z_0}(a,b)\right)\right| : (a,b)\in D_r(z_0)\times D_r(z_0)\right\}
\end{align*}
for every $(x,y)\in D_r(z_0)\times D_r(z_0)$. Therefore, if we can show that the series
\[
\sum_{j=1}^\infty \left[\frac{F_1}{\phi^{j+1}} + \sum_{k=0}^{j-1} \left|\frac{F_{j+1-k}}{\phi^{j+1}} -\frac{F_{j-k}}{\phi^j}\right|\right]
\]
converges (in $\R$), it will follow that our sequence of functions is uniformly Cauchy on $D_r(z_0)\times D_r(z_0)$.
\\

\noindent\underline{\textbf{Step 3:}} A well-known formula for the Fibonacci number $F_\ell$ is 
\[
F_\ell = \frac{\phi^\ell - (-1/\phi)^\ell}{\sqrt{5}} = \frac{\phi^\ell + (-1)^{\ell+1}\cdot \phi^{-\ell}}{\sqrt{5}}.
\]
Hence, for each $j\in \Z^+$ and each $0\leq k\leq j-1$, we find
\begin{align*}
\frac{F_{j+1-k}}{\phi^{j+1}} -\frac{F_{j-k}}{\phi^j} & = \frac1{\sqrt{5}}\cdot\left[\frac{\phi^{j+1-k}+(-1)^{j+2-k}\cdot \phi^{-(j+1-k)}}{\phi^{j+1}} -\frac{\phi^{j-k}+(-1)^{j+1-k}\cdot \phi^{-(j-k)}}{\phi^j}\right]\\
& = \frac1{\sqrt{5}}\cdot\left[\left(\phi^{-k}+(-1)^{j+2-k}\cdot \phi^{-2(j+1)}\cdot \phi^k\right) -\left(\phi^{-k}+(-1)^{j+1-k}\cdot \phi^{-2j}\cdot \phi^k\right)\right]\\
& = \frac{(-1)^{j+2-k}}{\sqrt{5}}\cdot \phi^{-2(j+1)}\left(1+\phi^2\right)\cdot \phi^k.
\end{align*}
The geometric series formula implies
\begin{align*}
\sum_{k=0}^{j-1} \left|\frac{F_{j+1-k}}{\phi^{j+1}} -\frac{F_{j-k}}{\phi^j}\right| & = \sum_{k=0}^{j-1} \left|\frac{(-1)^{j+2-k}}{\sqrt{5}}\cdot \phi^{-2(j+1)}\left(1+\phi^2\right)\cdot \phi^k\right| = \sum_{k=0}^{j-1}\frac{\phi^{-2(j+1)}\left(1+\phi^2\right)}{\sqrt{5}}\cdot \phi^k\\
& = \frac{\phi^{-2(j+1)}\left(1+\phi^2\right)}{\sqrt{5}}\cdot \frac{1-\phi^j}{1-\phi} = \frac1{\sqrt5}\cdot\frac{\phi^2+1}{\phi-1}\cdot \left(\frac1{\phi^{j+2}} - \frac1{\phi^{2j+2}}\right).
\end{align*}
The geometric series formula applies a second time to yield that our desired series converges, with
\begin{align*}
\sum_{j=1}^\infty \left[\frac{F_1}{\phi^{j+1}} + \sum_{k=0}^{j-1} \left|\frac{F_{j+1-k}}{\phi^{j+1}} -\frac{F_{j-k}}{\phi^j}\right|\right] & = \sum_{j=1}^\infty \left[\frac1{\phi^{j+1}} + \frac1{\sqrt5}\cdot\frac{\phi^2+1}{\phi-1}\cdot \left(\frac1{\phi^{j+2}} - \frac1{\phi^{2j+2}}\right) \right]\\
& = \frac1{\phi^2}\cdot \frac{1}{1-1/\phi} + \frac1{\sqrt5}\cdot\frac{\phi^2+1}{\phi-1}\cdot\left(\frac1{\phi^3}\cdot \frac{1}{1-1/\phi} - \frac1{\phi^4}\cdot \frac{1}{1-1/\phi^2}\right)\\
& = \frac{-5+\sqrt{5}+(5+\sqrt{5})\phi^2}{5\phi(\phi-1)^2(\phi+1)} = 2
\end{align*}
(the final simplification was obtained using the \texttt{Simplify} command in \texttt{Mathematica}).

Since the full series converges to $2$, its partial sums are a Cauchy sequence, i.e., the expression
\[
\sum_{j=n}^{m-1} \left[\frac{F_1}{\phi^{j+1}}+\sum_{k=0}^{j-1}\left|\frac{F_{j+1-k}}{\phi^{j+1}}-\frac{F_{j-k}}{\phi^j}\right|\right]
\]
can be made arbitrarily small as $m$ and $n$ become large. Therefore, by Step 2, the expression
\[
\left|\sum_{k=0}^{m-1} \frac{F_{m-k}}{\phi^m}\cdot \mc{L}\left(\mc{G}_{f,z_0}(\mc{S}_f^{\circ k}(x,y))\right) - \sum_{k=0}^{n-1} \frac{F_{n-k}}{\phi^n}\cdot \mc{L}\left(\mc{G}_{f,z_0}(\mc{S}_f^{\circ k}(x,y))\right)\right|
\]
can be made uniformly arbitrarily small for all $(x,y)\in D_r(z_0)\times D_r(z_0)$, as $m$ and $n$ become large. %So our original sequence of functions $g_n$ converges uniformly on all of $D_r(z_0)\times D_r(z_0)$.
Per our discussion in Step 2, this implies the $H_n$ converge uniformly on $D_r(z_0)\times D_r(z_0)$.

Define $\mc{H}_{f,z_0,w_0} : D_r(z_0)\times D_r(z_0) \to \C$ to be the limit
\[
\mc{H}_{f,z_0,w_0}(x,y) = \lim_{n\to\infty} H_n(x,y) = \lim_{n\to\infty}  \exp\left(\sum_{k=0}^{n-1} \frac{F_{n-k}}{\phi^n}\cdot \mc{L}\left(\mc{G}_{f,z_0}(\mc{S}_f^{\circ k}(x,y))\right)\right).
\]
Then $\mc{H}_{f,z_0,w_0}$ is holomorphic due to a theorem by Weierstrass, being a uniform limit of a sequence of holomorphic maps. It remains to verify that this $\mc{H}_{f,z_0,w_0}$ has the desired properties.
\\

\noindent\underline{\textbf{Step 4:}} First, we verify that
\[
\mc{H}_{f,z_0,w_0}(z_0,z_0) = \exp\left(\frac{5+3\sqrt{5}}{10}\cdot w_0\right) =: \left(\frac{f''(z_0)}{2f'(z_0)}\right)^{\frac{5+3\sqrt{5}}{10}}.
\]
Using the facts that $\mc{G}_{f,z_0}(z_0,z_0) = \exp(w_0)$ by definition of $w_0$, that $\mc{L}(\exp(w_0)) = w_0$ by definition of $\mc{L}$, and that $(z_0,z_0)$ is fixed under $\mc{S}_f$, we find
\begin{align*}
\mc{H}_{f,z_0,w_0}(z_0,z_0) & = \lim_{n\to \infty} H_n(z_0,z_0) = \lim_{n\to\infty}  \exp\left(\sum_{k=0}^{n-1} \frac{F_{n-k}}{\phi^n}\cdot \mc{L}\left(\mc{G}_{f,z_0}(\mc{S}_f^{\circ k}(z_0,z_0))\right)\right)\\
& = \lim_{n\to\infty} \exp\left(\sum_{k=0}^{n-1} \frac{F_{n-k}}{\phi^n}\cdot \mc{L}\left(\mc{G}_{f,z_0}(z_0,z_0)\right)\right) = \lim_{n\to\infty} \exp\left(\sum_{k=0}^{n-1} \frac{F_{n-k}}{\phi^n}\cdot \mc{L}\left(\exp(w_0)\right)\right)\\
& = \lim_{n\to\infty} \exp\left(\sum_{k=0}^{n-1} \frac{F_{n-k}}{\phi^n}\cdot w_0\right) = \exp\left(\left[\lim_{n\to\infty} \sum_{k=0}^{n-1} \frac{F_{n-k}}{\phi^n}\right]\cdot w_0\right).
\end{align*}
The last equality used continuity of $\exp$. By work similar to Step 3, we find
\begin{align*}
\sum_{k=0}^{n-1} \frac{F_{n-k}}{\phi^n} & = \sum_{k=0}^{n-1}\frac1{\sqrt5}\cdot\frac{\phi^{n-k}+(-1)^{n+1-k}\cdot \phi^{-(n-k)}}{\phi^n} = \frac1{\sqrt5}\cdot \sum_{k=0}^{n-1} \left(\phi^{-k} + (-1)^{n+1-k}\cdot \phi^{-2n}\cdot \phi^k\right)\\
& = \frac1{\sqrt{5}}\cdot \left(\frac{1-(1/\phi)^n}{1-1/\phi} + (-1)^{n+1}\cdot \phi^{-2n}\cdot \frac{1-(-1/\phi)^n}{1-(1/\phi)}\right).
\end{align*}
The limit is thus
\[
\lim_{n\to \infty} \sum_{k=0}^{n-1} \frac{F_{n-k}}{\phi^n} = \frac1{\sqrt{5}}\cdot \frac1{1-1/\phi} = \frac{5+3\sqrt{5}}{10}.
\]
So indeed,
\[
\mc{H}_{f,z_0,w_0}(z_0,z_0) = \exp\left(\left[\lim_{n\to\infty} \sum_{k=0}^{n-1} \frac{F_{n-k}}{\phi^n}\right]\cdot w_0\right) = \exp\left(\frac{5+3\sqrt{5}}{10}\cdot w_0\right),
\]
as desired.

Let $(x,y)\in D_r(z_0)\times D_r(z_0)$ be given. By Step 1 and the definition of $\mc{L}$, we have
\begin{align*}
\left|(\pi_1\circ\mc{S}_f^{\circ n})(x,y) - z_0\right| & = \left|\left(z_0+\mc{G}_n(x,y)(x-z_0)^{F_{n+1}}(y-z_0)^{F_n}\right)-z_0\right|\\
& = |\mc{G}_n(x,y)|\cdot |x-z_0|^{F_{n+1}}\cdot |y-z_0|^{F_n}\\
& = \left|\prod_{k=0}^{n-1} \left[\mc{G}_{f,z_0}(\mc{S}_f^{\circ k}(x,y))\right]^{F_{n-k}}\right|\cdot |x-z_0|^{F_{n+1}}\cdot |y-z_0|^{F_n}\\
& = \left|\prod_{k=0}^{n-1} \left[\exp\left(\mc{L}\left(\mc{G}_{f,z_0}(\mc{S}_f^{\circ k}(x,y))\right)\right)\right]^{F_{n-k}}\right|\cdot |x-z_0|^{F_{n+1}}\cdot |y-z_0|^{F_n}\\
& = \left|\exp\left(\sum_{k=0}^{n-1} F_{n-k}\cdot \mc{L}\left(\mc{G}_{f,z_0}(\mc{S}_f^{\circ k}(x,y))\right)\right)\right|\cdot |x-z_0|^{F_{n+1}}\cdot |y-z_0|^{F_n}\\
& = \exp\left(\re\left(\sum_{k=0}^{n-1} F_{n-k}\cdot \mc{L}\left(\mc{G}_{f,z_0}(\mc{S}_f^{\circ k}(x,y))\right)\right)\right)\cdot |x-z_0|^{F_{n+1}}\cdot |y-z_0|^{F_n}\\
& = \left[\exp\left(\re\left(\sum_{k=0}^{n-1} \frac{F_{n-k}}{\phi^n}\cdot \mc{L}\left(\mc{G}_{f,z_0}(\mc{S}_f^{\circ k}(x,y))\right)\right)\right)\right]^{\phi^n}\cdot |x-z_0|^{F_{n+1}}\cdot |y-z_0|^{F_n}\\
& = \left|\exp\left(\sum_{k=0}^{n-1} \frac{F_{n-k}}{\phi^n}\cdot \mc{L}\left(\mc{G}_{f,z_0}(\mc{S}_f^{\circ k}(x,y))\right)\right)\right|^{\phi^n} \cdot |x-z_0|^{F_{n+1}}\cdot |y-z_0|^{F_n}\\
& = \left|H_n(x,y)\right|^{\phi^n}\cdot |x-z_0|^{F_{n+1}}\cdot |y-z_0|^{F_n}.
\end{align*}
In particular, since
\[
\lim_{n\to \infty} \frac{F_{n+1}}{\phi^n} = \lim_{n\to\infty} \frac1{\sqrt5}\cdot \frac{\phi^{n+1} + (-1)^{n+1}\cdot\phi^{-(n+1)}}{\phi^n} = \lim_{n\to\infty}\frac1{\sqrt5}\cdot \left(\phi + \frac{(-1)^{n+1}}{\phi^{2n+1}}\right) = \frac\phi{\sqrt5}
\]
and
\[
\lim_{n\to \infty} \frac{F_n}{\phi^n} = \lim_{n\to\infty} \frac1{\sqrt5}\cdot \frac{\phi^n + (-1)^n\cdot\phi^{-n}}{\phi^n} = \lim_{n\to\infty}\frac1{\sqrt5}\cdot \left(1 + \frac{(-1)^n}{\phi^{2n}}\right) = \frac1{\sqrt5},
\]
it follows that 
\begin{align*}
\lim_{n\to \infty} \left|(\pi_1\circ\mc{S}_f^{\circ n})(x,y) - z_0\right|^{1/\phi^n} & = \lim_{n\to \infty}\left(\left|H_n(x,y)\right|^{\phi^n}\cdot |x-z_0|^{F_{n+1}}\cdot |y-z_0|^{F_n}\right)^{1/\phi^n}\\
& = \lim_{n\to \infty} |H_n(x,y)|\cdot |x-z_0|^{F_{n+1}/\phi^n}\cdot |y-z_0|^{F_n/\phi^n}\\
& = \left|\lim_{n\to\infty} H_n(x,y)\right| \cdot |x-z_0|^{\lim_{n\to\infty} F_{n+1}/\phi^n} \cdot |y-z_0|^{\lim_{n\to\infty} F_n/\phi^n}\\
& = |\mc{H}_{f,z_0,w_0}(x,y)|\cdot |x-z_0|^{\phi/\sqrt{5}}\cdot |y-z_0|^{1/\sqrt{5}},
\end{align*}
as desired. This completes item (1) of the theorem.
\\

\noindent\boxed{\textbf{(2).}} Suppose $V\subseteq \C^2$ is another open set satisfying $(z_0,z_0)\in V\subseteq \mc{A}(f,z_0)$ and $\mc{S}_f(V)\subseteq V$, and let $H : V\to \C$ be a holomorphic map such that
\[
H(z_0,z_0) = \exp\left(\frac{5+3\sqrt5}{10}\cdot w_0\right) = \mc{H}_{f,z_0,w_0}(z_0,z_0)
\]
and 
\[
\lim_{n\to \infty} \left|(\pi_1\circ\mc{S}_f^{\circ n})(x,y) - z_0\right|^{1/\phi^n} = |H(x,y)|\cdot |x-z_0|^{\phi/\sqrt{5}}\cdot |y-z_0|^{1/\sqrt{5}}
\]
for all $(x,y)\in V$. We show that 
\[
H(x,y) = \mc{H}_{f,z_0,w_0}(x,y)
\]
for all $(x,y)\in V\cap (D_r(z_0)\times D_r(z_0))$. 

Since $H$ and $\mc{H}_{f,z_0,w_0}$ are continuous at $(z_0,z_0)$, and since $V\cap (D_r(z_0)\times D_r(z_0))$ is open, we find there exists an $r' > 0$ such that
\[
D_{r'}(z_0)\times D_{r'}(z_0) \subseteq V\cap (D_r(z_0)\times D_r(z_0)) 
\]
and $(x,y)\in D_{r'}(z_0)\times D_{r'}(z_0)$ implies
\[
\left|H(x,y)-H(z_0,z_0)\right|, \left|\mc{H}_{f,z_0,w_0}(x,y) - \mc{H}_{f,z_0,w_0}(z_0,z_0)\right| < \frac12\cdot \left|\exp\left(\frac{5+3\sqrt{5}}{10}\cdot w_0\right)\right|.
\]
Since the disk $D_{\frac12\cdot\left|\exp\left(\frac{5+3\sqrt{5}}{10}\cdot w_0\right)\right|}\left(\exp\left(\frac{5+3\sqrt{5}}{10}\cdot w_0\right)\right)$ does not contain zero and is simply connected, there also exists a holomorphic branch $L : D_{\frac12\cdot\left|\exp\left(\frac{5+3\sqrt{5}}{10}\cdot w_0\right)\right|}\left(\exp\left(\frac{5+3\sqrt{5}}{10}\cdot w_0\right)\right) \to \C$ of the complex logarithm such that
\[
L\left(\exp\left(\frac{5+3\sqrt{5}}{10}\cdot w_0\right)\right) = \frac{5+3\sqrt{5}}{10}\cdot w_0.
\]
In particular, the compositions $L\circ H$ and $L\circ \mc{H}$ are both well-defined and holomorphic on $D_{r'}(z_0)\times D_{r'}(z_0)$. We will show that both compositions agree on this polydisk.

% \[
% H\left(D_{r'}(z_0)\times D_{r'}(z_0)\right), \mc{H}_{f,z_0,w_0}\left(D_{r'}(z_0)\times D_{r'}(z_0)\right) \subseteq D_{\frac12\cdot\left|\exp\left(\frac{5+3\sqrt{5}}{10}\cdot w_0\right)\right|}\left(\exp\left(\frac{5+3\sqrt{5}}{10}\cdot w_0\right)\right)
% \]

Observe that for all $(x,y)\in D_{r'}(z_0)\times D_{r'}(z_0)$, we have
\begin{align*}
|H(x,y)|\cdot |x-z_0|^{\phi/\sqrt{5}}\cdot |y-z_0|^{1/\sqrt{5}} & = \lim_{n\to \infty} \left|(\pi_1\circ\mc{S}_f^{\circ n})(x,y) - z_0\right|^{1/\phi^n} \\
& = |\mc{H}_{f,z_0,w_0}(x,y)|\cdot |x-z_0|^{\phi/\sqrt{5}}\cdot |y-z_0|^{1/\sqrt{5}}.
\end{align*}
Since $H$ and $\mc{H}_{f,z_0,w_0}$ are both continuous, we can take limits at points on the lines $\{z_0\}\times \C$ and $\C\times \{z_0\}$ to see this implies
\[
|H(x,y)| = |\mc{H}_{f,z_0,w_0}(x,y)|
\]
for every $(x,y)\in D_{r'}(z_0)\times D_{r'}(z_0)$. That is,
\begin{align*}
e^{\re(L(H(x,y)))} & = \left|\exp\left(L(H(x,y))\right)\right| = |H(x,y)| \\
& = |\mc{H}_{f,z_0,w_0}(x,y)| = \left|\exp\left(L(\mc{H}_{f,z_0,w_0}(x,y))\right)\right| =  e^{\re(L(\mc{H}_{f,z_0,w_0}(x,y)))}.
\end{align*}
Injectivity of $t\mapsto e^t$ on $\R$ thus gives 
\[
\re(L(H(x,y))) = \re(L(\mc{H}_{f,z_0,w_0}(x,y)))
\]
for all $(x,y)\in D_{r'}(z_0)\times D_{r'}(z_0)$. %It remains to compare the imaginary parts.

Now we show $H$ and $\mc{H}_{f,z_0,w_0}$ have the same power series about $(z_0,z_0)$ (note that both series must converge on all of $D_{r'}(z_0)\times D_{r'}(z_0)$ by holomorphicity). Write $x = \alpha+i\beta$ and $y = \gamma+i\delta$. %Then, using the Wirtinger derivatives and 
By the the Cauchy-Riemann equations, we find
% \begin{align*}
% \frac{\partial(L\circ H)}{\partial x} & = \frac12\left(\frac{\partial(L\circ H)}{\partial \alpha} - i\cdot\frac{\partial(L\circ H)}{\partial \beta}\right)\\
% & = 
% \end{align*}
\begin{align*}
\frac{\partial (L\circ H)}{\partial x} & = \frac{\partial(\re\circ L\circ H)}{\partial \alpha} + i\cdot \frac{\partial(\im\circ L\circ H)}{\partial \alpha} = \frac{\partial(\re\circ L\circ H)}{\partial \alpha} - i\cdot \frac{\partial(\re\circ L\circ H)}{\partial \beta}\\
& = \frac{\partial(\re\circ L\circ \mc{H}_{f,z_0,w_0})}{\partial \alpha} - i\cdot \frac{\partial(\re\circ L\circ \mc{H}_{f,z_0,w_0})}{\partial \beta} = \frac{\partial(\re\circ L\circ \mc{H}_{f,z_0,w_0})}{\partial \alpha} + i\cdot \frac{\partial(\im\circ L\circ \mc{H}_{f,z_0,w_0})}{\partial \alpha}\\
& = \frac{\partial (L\circ \mc{H}_{f,z_0,w_0})}{\partial x}
\end{align*}
and similarly
\begin{align*}
\frac{\partial (L\circ H)}{\partial y} & = \frac{\partial(\re\circ L\circ H)}{\partial \gamma} + i\cdot \frac{\partial(\im\circ L\circ H)}{\partial \gamma} = \frac{\partial(\re\circ L\circ H)}{\partial \gamma} - i\cdot \frac{\partial(\re\circ L\circ H)}{\partial \delta}\\
& = \frac{\partial(\re\circ L\circ \mc{H}_{f,z_0,w_0})}{\partial \gamma} - i\cdot \frac{\partial(\re\circ L\circ \mc{H}_{f,z_0,w_0})}{\partial \delta} = \frac{\partial(\re\circ L\circ \mc{H}_{f,z_0,w_0})}{\partial \gamma} + i\cdot \frac{\partial(\im\circ L\circ \mc{H}_{f,z_0,w_0})}{\partial \gamma}\\
& = \frac{\partial (L\circ \mc{H}_{f,z_0,w_0})}{\partial y}
\end{align*}
on all of $D_{r'}(z_0)\times D_{r'}(z_0)$, since the compositions $\re\circ L\circ H$ and $\re\circ L\circ \mc{H}_{f,z_0,w_0}$ agree there. From this point, we can take higher partial derivatives at $(z_0,z_0)$ and reveal
\[
\frac{\partial^d(L\circ H)}{\partial x^{d-k}\partial y^k}(z_0,z_0) = \frac{\partial^d(L\circ \mc{H}_{f,z_0,w_0})}{\partial x^{d-k}\partial y^k}(z_0,z_0)
\]
for all $d\in \Z^+$ and all $0\leq k\leq d$. Since $H(z_0,z_0) = \mc{H}_{f,z_0,w_0}(z_0,z_0)$, it follows that $H$ and $\mc{H}_{f,z_0,w_0}$ have the same power series about $(z_0,z_0)$. Therefore,
\[
H(x,y) = \mc{H}_{f,z_0,w_0}(x,y)
\]
for all $(x,y)\in D_{r'}(z_0)\times D_{r'}(z_0)$, on comparing these series.

The Uniqueness Principle for holomorphic maps now applies to let us conclude
\[
H(x,y) = \mc{H}_{f,z_0,w_0}(x,y)
\]
for all $(x,y) \in V\cap (D_r(z_0)\times D_r(z_0))$, as desired.
\\

\noindent\boxed{\textbf{(3) and (4).}} \underline{\textbf{Step 1:}} Choose any $w_0\in \log\left(\frac{f''(z_0)}{2f'(z_0)}\right)$, and let $V\subseteq \C^2$ and $\mc{H}_{f,z_0,w_0} : V\to \C$ be any pair of domain and function satisfying (1). Then $(z_0,z_0)\in V\subseteq \mc{A}(f,z_0)$, so a general property of attracting basins implies we have
\[
\mc{A}(f,z_0) = \bigcup_{n=0}^\infty \mc{S}_f^{\circ(-n)}(V).
\]
We thus obtain an index function $\mc{N} : \mc{A}(f,z_0) \to \Z_{\geq0}$ given by
\[
\mc{N}(x,y) = \min\left\{n\in \Z_{\geq0} : (x,y) \in \mc{S}_f^{\circ (-n)}(V)\right\}.
\] 
Note that $\mc{N}^{-1}(\{0\}) = V$.

Using induction on the desired functional equation in (4), we define $\widehat{\mc{H}}_{f,z_0} : \mc{A}(f,z_0) \to \R$ by 
\[
\widehat{\mc{H}}_{f,z_0}(x,y) = \left|\mc{H}_{f,z_0,w_0}\left(\mc{S}_f^{\circ \mc{N}(x,y)}(x,y)\right)\right|^{\phi^{-\mc{N}(x,y)}} \cdot \prod_{k=0}^{\mc{N}(x,y)-1} \left|\mc{G}_{f,z_0}(\mc{S}_f^{\circ k}(x,y))\right|^{\phi^{-k}/\sqrt{5}}
\]
(adopting the convention that the product is just $1$ if $\mc{N}(x,y) = 0$). It is well-defined since $\mc{S}_f^{\circ \mc{N}(x,y)}(x,y) \in V$ by definition of $\mc{N}$. We show that this $\widehat{\mc{H}}_{f,z_0}$ is continuous on $\mc{A}(f,z_0)$, real-analytic on the preimage $\widehat{\mc{H}}_{f,z_0}^{-1}(\R\setminus\{0\})$, and it satisfies the functional equation from (4) as well as the equation
\[
\lim_{n\to \infty}\left|(\pi_1\circ\mc{S}_f^{\circ n})(x,y)-z_0\right|^{1/\phi^n} = \widehat{\mc{H}}_{f,z_0}(x,y)\cdot |x-z_0|^{\phi/\sqrt{5}}\cdot |y-z_0|^{1/\sqrt5}
\]
for all $(x,y)\in \mc{A}(f,z_0)$ (in particular, the LHS limit always exists).
\\

\noindent\underline{\textbf{Step 2:}} For continuity, we show via induction on $n\in\Z_{\geq0}$ that we have
\[
\widehat{\mc{H}}_{f,z_0}(x,y) = \left|\mc{H}_{f,z_0,w_0}\left(\mc{S}_f^{\circ n}(x,y)\right)\right|^{\phi^{-n}} \cdot \prod_{k=0}^{n-1} \left|\mc{G}_{f,z_0}(\mc{S}_f^{\circ k}(x,y))\right|^{\phi^{-k}/\sqrt{5}}
\]
for all $(x,y)\in \mc{S}_f^{\circ(-n)}(V)$ (this is slightly different from the definition of $\widehat{\mc{H}}_{f,z_0}$, with $\mc{N}(x,y)$ replaced by $n$). The base case $n = 0$ is simply the fact that $\mc{N}^{-1}(\{0\}) = V$, hence $\mc{N}(x,y) = 0$ for any $(x,y)\in V = \mc{S}_f^{\circ(-0)}(V)$. For the inductive step, suppose the claim holds when $n = m$, for some $m\in \Z_{\geq0}$. 

Let $(x,y)\in \mc{S}_f^{\circ(-(m+1))}(V)$ be given. Suppose first that $(x,y)\notin \mc{S}_f^{\circ(-m)}(V)$. Then, since $\mc{S}_f(V)\subseteq V$, we must have $\mc{N}(x,y) > m$ (lest otherwise we have $\mc{S}_f^{\circ m}(x,y)\in V$, contradiction). But we also have $\mc{N}(x,y) \leq m+1$ by minimality. Therefore, $\mc{N}(x,y) = m+1$, and $\widehat{\mc{H}}_{f,z_0}(x,y)$ has the form
\[
\widehat{\mc{H}}_{f,z_0}(x,y) = \left|\mc{H}_{f,z_0,w_0}\left(\mc{S}_f^{\circ (m+1)}(x,y)\right)\right|^{\phi^{-(m+1)}} \cdot \prod_{k=0}^m \left|\mc{G}_{f,z_0}(\mc{S}_f^{\circ k}(x,y))\right|^{\phi^{-k}/\sqrt{5}},
\]
as desired. Suppose instead that $(x,y)\in \mc{S}_f^{\circ(-m)}(V)$. Then by the inductive hypothesis, we have
\[
\widehat{\mc{H}}_{f,z_0}(x,y) = \left|\mc{H}_{f,z_0,w_0}\left(\mc{S}_f^{\circ m}(x,y)\right)\right|^{\phi^{-m}} \cdot \prod_{k=0}^{m-1} \left|\mc{G}_{f,z_0}(\mc{S}_f^{\circ k}(x,y))\right|^{\phi^{-k}/\sqrt{5}}.
\]
Observe that by (1), we have
\begin{align*}
\left|\mc{H}_{f,z_0,w_0}(a,b)\right| & \cdot |a-z_0|^{\phi/\sqrt5}\cdot |b-z_0|^{1/\sqrt5}\\
& = \lim_{n\to\infty} \left|(\pi_1\circ \mc{S}_f^{\circ n})(a,b) - z_0\right|^{1/\phi^n} = \lim_{n\to\infty} \left|(\pi_1\circ \mc{S}_f^{\circ (n+1)})(a,b) - z_0\right|^{1/\phi^{n+1}}\\
& = \left(\lim_{n\to\infty} \left|(\pi_1\circ \mc{S}_f^{\circ n})(\mc{S}(a,b)) - z_0\right|^{1/\phi^n}\right)^{1/\phi}\\
& = \left[\left|\mc{H}_{f,z_0,w_0}(\mc{S}_f(a,b))\right|\cdot \left|\pi_1(\mc{S}_f(a,b))-z_0\right|^{\phi/\sqrt5}\cdot \left|\pi_2(\mc{S}_f(a,b))-z_0\right|^{1/\sqrt5}\right]^{1/\phi}\\
& = \left[\left|\mc{H}_{f,z_0,w_0}(\mc{S}_f(a,b))\right|\cdot \left|\mc{G}_{f,z_0}(a,b)(a-z_0)(b-z_0)\right|^{\phi/\sqrt5}\cdot |a-z_0|^{1/\sqrt5}\right]^{1/\phi}\\
& = \left|\mc{H}_{f,z_0,w_0}(\mc{S}_f(a,b))\right|^{1/\phi} \cdot |\mc{G}_{f,z_0}(a,b)|^{1/\sqrt5}\cdot |a-z_0|^{(1+1/\phi)/\sqrt5}\cdot |b-z_0|^{1/\sqrt5}\\
& = \left|\mc{H}_{f,z_0,w_0}(\mc{S}_f(a,b))\right|^{1/\phi} \cdot |\mc{G}_{f,z_0}(a,b)|^{1/\sqrt5}\cdot |a-z_0|^{\phi/\sqrt5}\cdot |b-z_0|^{1/\sqrt5}
\end{align*}
for all $(a,b)\in V$. By continuity of $\mc{H}_{f,z_0,w_0}$ and of $\mc{G}_{f,z_0}$, we can take limits at points on the lines $\{z_0\}\times \C$ and $\C\times \{z_0\}$ to obtain that
\[
|\mc{H}_{f,z_0,w_0}(a,b)| = \left|\mc{H}_{f,z_0,w_0}(\mc{S}_f(a,b))\right|^{1/\phi} \cdot |\mc{G}_{f,z_0}(a,b)|^{1/\sqrt5}
\]
for all $(a,b)\in V$. In our case, since $\mc{S}_f^{\circ m}(x,y)\in V$, we obtain
\begin{align*}
\widehat{\mc{H}}_{f,z_0}(x,y) & = \left|\mc{H}_{f,z_0,w_0}\left(\mc{S}_f^{\circ m}(x,y)\right)\right|^{\phi^{-m}} \cdot \prod_{k=0}^{m-1} \left|\mc{G}_{f,z_0}(\mc{S}_f^{\circ k}(x,y))\right|^{\phi^{-k}/\sqrt{5}}\\
& = \left[\left|\mc{H}_{f,z_0,w_0}\left(\mc{S}_f^{\circ (m+1)}(x,y)\right)\right|^{1/\phi} \cdot \left|\mc{G}_{f,z_0}\left(\mc{S}_f^{\circ m}(x,y)\right)\right|^{1/\sqrt5}\right]^{\phi^{-m}} \cdot \prod_{k=0}^{m-1} \left|\mc{G}_{f,z_0}(\mc{S}_f^{\circ k}(x,y))\right|^{\phi^{-k}/\sqrt{5}}\\
& = \left|\mc{H}_{f,z_0,w_0}\left(\mc{S}_f^{\circ (m+1)}(x,y)\right)\right|^{1/\phi^{m+1}} \cdot \left|\mc{G}_{f,z_0}\left(\mc{S}_f^{\circ m}(x,y)\right)\right|^{\phi^{-m}/\sqrt5} \cdot \prod_{k=0}^{m-1} \left|\mc{G}_{f,z_0}(\mc{S}_f^{\circ k}(x,y))\right|^{\phi^{-k}/\sqrt{5}}\\
& = \left|\mc{H}_{f,z_0,w_0}\left(\mc{S}_f^{\circ (m+1)}(x,y)\right)\right|^{\phi^{-(m+1)}} \cdot \prod_{k=0}^m \left|\mc{G}_{f,z_0}(\mc{S}_f^{\circ k}(x,y))\right|^{\phi^{-k}/\sqrt{5}},
\end{align*}
as desired.

By induction, the claim holds for all $n\in \Z_{\geq0}$. Since each preimage $\mc{S}_f^{\circ (-n)}(V)$ is open, and
\[
\widehat{\mc{H}}_{f,z_0}(x,y) = \left|\mc{H}_{f,z_0,w_0}\left(\mc{S}_f^{\circ n}(x,y)\right)\right|^{\phi^{-n}} \cdot \prod_{k=0}^{n-1} \left|\mc{G}_{f,z_0}(\mc{S}_f^{\circ k}(x,y))\right|^{\phi^{-k}/\sqrt{5}}
\]
for all $(x,y)\in \mc{S}_f^{\circ (-n)}(V)$, it follows from continuity of $\mc{H}_{f,z_0,w_0}$, of $\mc{G}_{f,z_0}$, and of $\mc{S}_f$ that $\widehat{\mc{H}}_{f,z_0}$ is continuous on $\mc{S}_f^{\circ (-n)}(V)$. Therefore, $\widehat{\mc{H}}_{f,z_0}$ is continuous on all of $\mc{A}(f,z_0)$, as claimed. Note that this implies the preimage $\widehat{\mc{H}}_{f,z_0}^{-1}(\R\setminus\{0\}) \subseteq \mc{A}(f,z_0)$ is open in $\C^2$.

For real-analyticity on $\widehat{\mc{H}}_{f,z_0}^{-1}(\R\setminus\{0\})$, note that it is sufficient to prove real-analyticity on each open subset $\mc{S}_f^{\circ (-n)}(V) \cap \widehat{\mc{H}}_{f,z_0}^{-1}(\R\setminus\{0\})$, for $n\in \Z_{\geq0}$. This is immediate from our above formula for $\widehat{\mc{H}}_{f,z_0}$ on $\mc{S}_f^{\circ (-n)}(V)$, due to the real-analyticity of the maps $(x,y)\mapsto \left|\mc{H}_{f,z_0,w_0}\left(\mc{S}_f^{\circ n}(x,y)\right)\right|^{\phi^{-n}}$ and $(x,y) \mapsto \left|\mc{G}_{f,z_0}(\mc{S}_f^{\circ k}(x,y))\right|^{\phi^{-k}/\sqrt{5}}$ away from the points where they are zero. So indeed, $\widehat{\mc{H}}_{f,z_0}$ is real-analytic on $\widehat{\mc{H}}_{f,z_0}^{-1}(\R\setminus\{0\})$.
\\

\noindent\underline{\textbf{Step 3:}} For the limit, let $(x,y)\in \mc{A}(f,z_0)$ be given. If $(x,y) \in V$, then there is nothing to show, because $\widehat{\mc{H}}_{f,z_0}(x,y) = |\mc{H}_{f,z_0,w_0}(x,y)|$, and we can apply (1). Suppose instead that $(x,y)\notin V$. Then $\mc{N}_0 := \mc{N}(x,y) > 0$. Because $\mc{S}_f^{\circ \mc{N}_0}(x,y)\in V$, we find via (1) that
\begin{align*}
\lim_{n\to\infty} & \left|(\pi_1\circ \mc{S}_f^{\circ n})\left(\mc{S}_f^{\circ \mc{N}_0}(x,y)\right)-z_0\right|^{1/\phi^n} \\
& = \left|\mc{H}_{f,z_0,w_0}\left(\mc{S}_f^{\circ \mc{N}_0}(x,y)\right)\right| \cdot \left|\pi_1\left(\mc{S}_f^{\circ \mc{N}_0}(x,y)\right)-z_0\right|^{\phi/\sqrt{5}}\cdot \left|\pi_2\left(\mc{S}_f^{\circ \mc{N}_0}(x,y)\right)-z_0\right|^{1/\sqrt5}\\
& = \left|\mc{H}_{f,z_0,w_0}\left(\mc{S}_f^{\circ \mc{N}_0}(x,y)\right)\right| \cdot \left|\mc{G}_{\mc{N}_0}(x,y)(x-z_0)^{F_{\mc{N}_0+1}}(y-z_0)^{F_{\mc{N}_0}}\right|^{\phi/\sqrt5}\\
& \ \ \ \ \ 
\cdot \left|\mc{G}_{\mc{N}_0-1}(x,y)(x-z_0)^{F_{\mc{N}_0}}(y-z_0)^{F_{\mc{N}_0-1}}\right|^{1/\sqrt5}\\
& = \left|\mc{H}_{f,z_0,w_0}\left(\mc{S}_f^{\circ \mc{N}_0}(x,y)\right)\right| \cdot \left|\prod_{k=0}^{\mc{N}_0-1} \left[\mc{G}_{f,z_0}(\mc{S}_f^{\circ k}(x,y))\right]^{F_{\mc{N}_0-k}}\right|^{\phi/\sqrt5} \cdot \left|\prod_{k=0}^{\mc{N}_0-2} \left[\mc{G}_{f,z_0}(\mc{S}_f^{\circ k}(x,y))\right]^{F_{\mc{N}_0-1-k}}\right|^{1/\sqrt5}\\
& \ \ \ \ \ \cdot |x-z_0|^{(1/\sqrt5)(\phi\cdot F_{\mc{N}_0+1}+F_{\mc{N}_0})}\cdot |y-z_0|^{(1/\sqrt5)(\phi\cdot F_{\mc{N}_0}+F_{\mc{N}_0-1})}\\
& = \left|\mc{H}_{f,z_0,w_0}\left(\mc{S}_f^{\circ \mc{N}_0}(x,y)\right)\right| \cdot\left|\mc{G}_{f,z_0}\left(\mc{S}_f^{\circ(\mc{N}_0-1)}(x,y)\right)\right|^{\phi\cdot F_1/\sqrt5} \\
& \ \ \ \ \ \cdot \left(\prod_{k=0}^{\mc{N}_0-2} \left|\mc{G}_{f,z_0}(\mc{S}_f^{\circ k}(x,y))\right|^{(1/\sqrt5)(\phi\cdot F_{\mc{N}_0-k}+F_{\mc{N}_0-1-k})}\right) \\
& \ \ \ \ \ \cdot |x-z_0|^{(1/\sqrt5)(\phi\cdot F_{\mc{N}_0+1}+F_{\mc{N}_0})}\cdot |y-z_0|^{(1/\sqrt5)(\phi\cdot F_{\mc{N}_0}+F_{\mc{N}_0-1})}\\
& = \left|\mc{H}_{f,z_0,w_0}\left(\mc{S}_f^{\circ \mc{N}_0}(x,y)\right)\right| \cdot \left(\prod_{k=0}^{\mc{N}_0-1} \left|\mc{G}_{f,z_0}(\mc{S}_f^{\circ k}(x,y))\right|^{(1/\sqrt5)(\phi\cdot F_{\mc{N}_0-k}+F_{\mc{N}_0-1-k})}\right)\\
& \ \ \ \ \ \cdot |x-z_0|^{(1/\sqrt5)(\phi\cdot F_{\mc{N}_0+1}+F_{\mc{N}_0})}\cdot |y-z_0|^{(1/\sqrt5)(\phi\cdot F_{\mc{N}_0}+F_{\mc{N}_0-1})}.
\end{align*}
Note that for each $j\in \Z_{\geq0}$, we have
\begin{align*}
\frac1{\sqrt5}(\phi\cdot F_{j+1}+F_j) & = \frac1{\sqrt5}\left(\phi\cdot \frac{\phi^{j+1}+(-1)^{j+2}\cdot \phi^{-(j+1)}}{\sqrt5} + \frac{\phi^j + (-1)^{j+1}\cdot \phi^{-j}}{\sqrt5}\right)\\
& = \frac15\left(\phi^{j+2} + (-1)^{j+2}\cdot \phi^{-j} + \phi^j + (-1)^{j+1}\cdot \phi^{-j}\right) = \frac15\left(\phi^{j+2}+\phi^j\right)\\
& = \frac{\phi^{j+1}}5\left(\phi+\frac1\phi\right) = \frac{\phi^{j+1}}5\left(\frac{1+\sqrt5}{2}-\frac{1-\sqrt5}2\right)\\
& = \frac{\phi^{j+1}}{\sqrt5},
\end{align*}
whence
\begin{align*}
\lim_{n\to\infty} & \left|(\pi_1\circ \mc{S}_f^{\circ n})\left(\mc{S}_f^{\circ \mc{N}_0}(x,y)\right)-z_0\right|^{1/\phi^n} \\
& = \left|\mc{H}_{f,z_0,w_0}\left(\mc{S}_f^{\circ \mc{N}_0}(x,y)\right)\right| \cdot \left(\prod_{k=0}^{\mc{N}_0-1} \left|\mc{G}_{f,z_0}(\mc{S}_f^{\circ k}(x,y))\right|^{(1/\sqrt5)(\phi\cdot F_{\mc{N}_0-k}+F_{\mc{N}_0-1-k})}\right)\\
& \ \ \ \ \ \cdot |x-z_0|^{(1/\sqrt5)(\phi\cdot F_{\mc{N}_0+1}+F_{\mc{N}_0})}\cdot |y-z_0|^{(1/\sqrt5)(\phi\cdot F_{\mc{N}_0}+F_{\mc{N}_0-1})}\\
& = \left|\mc{H}_{f,z_0,w_0}\left(\mc{S}_f^{\circ \mc{N}_0}(x,y)\right)\right| \cdot \left(\prod_{k=0}^{\mc{N}_0-1} \left|\mc{G}_{f,z_0}(\mc{S}_f^{\circ k}(x,y))\right|^{\phi^{\mc{N}_0-k}/\sqrt5}\right)\cdot |x-z_0|^{\phi^{\mc{N}_0+1}/\sqrt5}\cdot |y-z_0|^{\phi^{\mc{N}_0}/\sqrt5}\\
& = \left[\left|\mc{H}_{f,z_0,w_0}\left(\mc{S}_f^{\circ \mc{N}_0}(x,y)\right)\right|^{\phi^{-\mc{N}_0}} \cdot \left(\prod_{k=0}^{\mc{N}_0-1} \left|\mc{G}_{f,z_0}(\mc{S}_f^{\circ k}(x,y))\right|^{\phi^{-k}/\sqrt5}\right)\cdot |x-z_0|^{\phi/\sqrt5}\cdot |y-z_0|^{1/\sqrt5}\right]^{\phi^{\mc{N}_0}}\\
& = \left[\widehat{\mc{H}}_{f,z_0}(x,y)\cdot |x-z_0|^{\phi/\sqrt5}\cdot |y-z_0|^{1/\sqrt5}\right]^{\phi^{\mc{N}_0}}.
\end{align*}
Raising this equation to the $\phi^{-\mc{N}_0}$ power therefore yields
\begin{align*}
\widehat{\mc{H}}_{f,z_0}(x,y)\cdot |x-z_0|^{\phi/\sqrt5}\cdot |y-z_0|^{1/\sqrt5} & = \left(\lim_{n\to\infty} \left|(\pi_1\circ \mc{S}_f^{\circ n})\left(\mc{S}_f^{\circ \mc{N}_0}(x,y)\right)-z_0\right|^{1/\phi^n}\right)^{\phi^{-\mc{N}_0}}\\
& = \lim_{n\to \infty} \left|(\pi_1\circ \mc{S}_f^{\circ (n+\mc{N}_0)})(x,y)-z_0\right|^{1/\phi^{n+\mc{N}_0}}\\
& = \lim_{n\to \infty} \left|(\pi_1\circ \mc{S}_f^{\circ n})(x,y)-z_0\right|^{1/\phi^n},
\end{align*}
as desired.
\\

\noindent\underline{\textbf{Step 4:}} Finally, for the functional equation from (4), we use the limit of Step 3 to adapt the work of Step 2 and find
\begin{align*}
\widehat{\mc{H}}_{f,z_0}(x,y) & \cdot |x-z_0|^{\phi/\sqrt5}\cdot |y-z_0|^{1/\sqrt5}\\
& = \lim_{n\to\infty} \left|(\pi_1\circ \mc{S}_f^{\circ n})(x,y) - z_0\right|^{1/\phi^n} = \lim_{n\to\infty} \left|(\pi_1\circ \mc{S}_f^{\circ (n+1)})(x,y) - z_0\right|^{1/\phi^{n+1}}\\
& = \left(\lim_{n\to\infty} \left|(\pi_1\circ \mc{S}_f^{\circ n})(\mc{S}(x,y)) - z_0\right|^{1/\phi^n}\right)^{1/\phi}\\
& = \left[\widehat{\mc{H}}_{f,z_0}(\mc{S}_f(x,y))\cdot \left|\pi_1(\mc{S}_f(x,y))-z_0\right|^{\phi/\sqrt5}\cdot \left|\pi_2(\mc{S}_f(x,y))-z_0\right|^{1/\sqrt5}\right]^{1/\phi}\\
& = \left[\widehat{\mc{H}}_{f,z_0}(\mc{S}_f(x,y))\cdot \left|\mc{G}_{f,z_0}(x,y)(x-z_0)(y-z_0)\right|^{\phi/\sqrt5}\cdot |x-z_0|^{1/\sqrt5}\right]^{1/\phi}\\
& = \left[\widehat{\mc{H}}_{f,z_0}(\mc{S}_f(x,y))\right]^{1/\phi} \cdot |\mc{G}_{f,z_0}(x,y)|^{1/\sqrt5}\cdot |x-z_0|^{(1+1/\phi)/\sqrt5}\cdot |y-z_0|^{1/\sqrt5}\\
& = \left[\widehat{\mc{H}}_{f,z_0}(\mc{S}_f(x,y))\right]^{1/\phi} \cdot |\mc{G}_{f,z_0}(x,y)|^{1/\sqrt5}\cdot |x-z_0|^{\phi/\sqrt5}\cdot |y-z_0|^{1/\sqrt5}
\end{align*}
for all $(x,y)\in \mc{A}(f,z_0)$. By continuity of $\widehat{\mc{H}}_{f,z_0}$ and of $\mc{G}_{f,z_0}$, we can take limits at points on the lines $\{z_0\}\times \C$ and $\C\times \{z_0\}$ to obtain that
\[
\widehat{\mc{H}}_{f,z_0}(x,y) = \left[\widehat{\mc{H}}_{f,z_0}(\mc{S}_f(x,y))\right]^{1/\phi} \cdot |\mc{G}_{f,z_0}(x,y)|^{1/\sqrt5}
\]
for all $(x,y)\in \mc{A}(f,z_0)$, as claimed. Combining the work of all the above steps, it follows that $\widehat{\mc{H}}_{f,z_0}$ has all the properties we want for (3) and (4); it remains to prove uniqueness.
\\

\noindent\underline{\textbf{Step 5:}} Let $\widehat{H} : \mc{A}(f,z_0) \to \R$ be another continuous function such that
\[
\lim_{n\to \infty} \left|(\pi_1\circ\mc{S}_f^{\circ n})(x,y)-z_0\right|^{1/\phi^n} = \widehat{H}(x,y)\cdot |x-z_0|^{\phi/\sqrt5}\cdot |y-z_0|^{1/\sqrt5}
\]
for all $(x,y)\in \mc{A}(f,z_0)$. Then certainly
\[
\widehat{H}(x,y)\cdot |x-z_0|^{\phi/\sqrt5}\cdot |y-z_0|^{1/\sqrt5} = \widehat{\mc{H}}_{f,z_0}(x,y)\cdot |x-z_0|^{\phi/\sqrt5}\cdot |y-z_0|^{1/\sqrt5}
\]
by Step 3. Once again, by continuity of $\widehat{H}$ and $\widehat{\mc{H}}_{f,z_0}$, we can take limits at points on the lines $\{z_0\}\times \C$ and $\C\times \{z_0\}$ to obtain that
\[
\widehat{H}(x,y) = \widehat{\mc{H}}_{f,z_0}(x,y).
\]
So $\widehat{H} = \widehat{\mc{H}}_{f,z_0}$, and our function is unique.
\end{proof}

This completes most of the hard work needed to show the potential $h_{f,z_0}$ is well-defined. However, we are not done yet: the limit $\lim_{n\to\infty} \left|(\pi_1\circ \mc{S}_f^{\circ n})(x,y)-z_0\right|^{1/\phi^n}$ considered in Theorem \hyperlink{3:BöttcherFunction}{3.1} does not, a priori, equal a $\lim_{n\to\infty} \left\|\mc{S}_f^{\circ n}(x,y)-(z_0,z_0)\right\|^{1/\phi^n}$ for any norm $\|\cdot\| : \C^2\to \C$, because the map $(x,y) \mapsto |x|$ is not itself a norm. 

We are fortunate, however, that we know the second coordinate $(\pi_2\circ \mc{S}_f^{\circ n})(x,y)$ in terms of the first due to Proposition \hyperlink{2:Iterates}{2.3}. Thus, we can compute the limit $\lim_{n\to\infty} \left|(\pi_2\circ \mc{S}_f^{\circ n})(x,y)-z_0\right|^{1/\phi^n}$ and combine our results in Theorem \hyperlink{3:BöttcherFunction}{3.1} to arrive at the true value of the potential.

\begin{center}
\begin{minipage}{37em}
\noindent\textbf{Corollary \hyperlink{3:Potential}{3.2} (Potential Function)} \textit{Suppose we have $f''(z_0) \neq 0$. Let $\widehat{\mc{H}}_{f,z_0} : \mc{A}(f,z_0) \to \R$ be the Böttcher-type modulus about $(z_0,z_0)$, from Theorem \hyperlink{3:BöttcherFunction}{3.1}.}
\begin{itemize}
    \item[\textbf{\textit{(1).}}] \textit{Define a map $h_{f,z_0} : \mc{A}(f,z_0)\to \R$ by
    \begin{align*}
    h_{f,z_0}(x,y) & = \max\left\{\widehat{\mc{H}}_{f,z_0}(x,y)\cdot |x-z_0|^{\phi/\sqrt5}\cdot |y-z_0|^{1/\sqrt5},\right.\\
    & \ \ \ \ \ \ \ \ \ \ \ \left.\left[\widehat{\mc{H}}_{f,z_0}(x,y)\right]^{1/\phi}\cdot |x-z_0|^{1/\sqrt5}\cdot |y-z_0|^{\phi^{-1}/\sqrt5}\right\}.
    \end{align*}
    Then $h_{f,z_0}$ is continuous and satisfies
    \[
    h_{f,z_0}(\mc{S}_f(x,y)) = [h_{f,z_0}(x,y)]^\phi
    \]
    for all $(x,y)\in \mc{A}(f,z_0)$.
    \item[\textbf{(2).}] For every norm $\|\cdot\| : \C^2\to \R$ on $\C^2$, we have 
    \[
    h_{f,z_0}(x,y) = \lim_{n\to \infty} \left\|\mc{S}_f^{\circ n}(x,y) - (z_0,z_0)\right\|^{1/\phi^n}
    \]
    for all $(x,y)\in \mc{A}(f,z_0)$. Therefore, $h_{f,z_0}$ is the potential function of $\mc{S}_f$ near $(z_0,z_0)$.
    \item[\textbf{(3).}] We have
    \[
    h_{f,z_0}^{-1}(\{0\}) = \bigcup_{n=0}^\infty \mc{S}_f^{\circ(-n)}(\{(z_0,z_0)\}),
    \]
    that is, the zero equipotential is equal to the set of iterated preimages of $(z_0,z_0)$.
    \item[\textbf{(4).}] More generally, we always have $0\leq h_{f,z_0}(x,y) < 1$, hence
    \[
    h_{f,z_0}(x,y) = \left[\widehat{\mc{H}}_{f,z_0}(x,y)\right]^{1/\phi}\cdot |x-z_0|^{1/\sqrt5}\cdot |y-z_0|^{\phi^{-1}/\sqrt5}
    \]
    for all $(x,y)\in \mc{A}(f,z_0)$.}
\end{itemize}
\end{minipage}
\end{center}
\begin{proof}
\boxed{\textbf{(1).}} %That $h_{f,z_0}$ is continuous follows from continuity of $\widehat{\mc{H}}_{f,z_0}$. 
Consider the temporary function $\widetilde{h} : \mc{A}(f,z_0) \to \R$ given by
\[
\widetilde{h}(x,y) = \widehat{\mc{H}}_{f,z_0}(x,y)\cdot |x-z_0|^{\phi/\sqrt5}\cdot |y-z_0|^{1/\sqrt5}.
\]
Then we see
\[
h_{f,z_0}(x,y) = \max\left\{\widetilde{h}(x,y), \left[\widetilde{h}(x,y)\right]^{1/\phi}\right\} = \begin{cases}
    \left[\widetilde{h}(x,y)\right]^{1/\phi} & \text{if } (x,y)\in \widetilde{h}^{-1}\left([0,1]\right)\\
    \widetilde{h}(x,y) & \text{if } (x,y)\in \widetilde{h}^{-1}\left([1,\infty)\right)
\end{cases}
\]
since $0 < 1/\phi < 1$ and $\widetilde{h}$ is always nonnegative. By continuity of $\widehat{\mc{H}}_{f,z_0}$, we see $\widetilde{h}$ is continuous, hence so is $h_{f,z_0}$ by pasting along the overlap region $\widetilde{h}^{-1}(\{1\})$ (whereupon $h_{f,z_0}$ is identically $1$). 

For the functional equation, we use the functional equation for $\widehat{\mc{H}}_{f,z_0}$: given $(x,y)\in \mc{A}(f,z_0)$, we have
\begin{align*}
\widetilde{h}(x,y) & = \widehat{\mc{H}}_{f,z_0}(x,y)\cdot |x-z_0|^{\phi/\sqrt5}\cdot |y-z_0|^{1/\sqrt5}\\
& = \left[\widehat{\mc{H}}_{f,z_0}(\mc{S}_f(x,y))\right]^{1/\phi}\cdot |\mc{G}_{f,z_0}(x,y)|^{1/\sqrt5}\cdot |x-z_0|^{\phi/\sqrt5}\cdot |y-z_0|^{1/\sqrt5}\\
& = \left[\widehat{\mc{H}}_{f,z_0}(\mc{S}_f(x,y))\right]^{1/\phi}\cdot \left|\mc{G}_{f,z_0}(x,y)(x-z_0)(y-z_0)\right|^{1/\sqrt5}\cdot |x-z_0|^{(1/\sqrt5)(\phi-1)}\\
& = \left[\widehat{\mc{H}}_{f,z_0}(\mc{S}_f(x,y))\cdot \left|\mc{G}_{f,z_0}(x,y)(x-z_0)(y-z_0)\right|^{\phi/\sqrt5}\cdot |x-z_0|^{(1/\sqrt5)(\phi^2-\phi)}\right]^{1/\phi}\\
& = \left[\widehat{\mc{H}}_{f,z_0}(\mc{S}_f(x,y))\cdot \left|\pi_1\left(\mc{S}_f(x,y)\right)-z_0\right|^{\phi/\sqrt5}\cdot \left|\pi_2\left(\mc{S}_f(x,y)\right)-z_0\right|^{1/\sqrt5}\right]^{1/\phi}\\
& = \left[\widetilde{h}(\mc{S}_f(x,y))\right]^{1/\phi}
\end{align*}
(where $\pi_1,\pi_2 : \C^2\to \C$ are the coordinate projections). Raising to the power $\phi$, we obtain
\[
\widetilde{h}(\mc{S}_f(x,y)) = \left[\widetilde{h}(x,y)\right]^\phi.
\]
In particular, we see that the three sets $\widetilde{h}^{-1}([0,1))$, $\widetilde{h}^{-1}(\{1\})$, and $\widehat{h}^{-1}([1,\infty))$ are invariant under $\mc{S}_f$, so we obtain
\begin{align*}
h_{f,z_0}(\mc{S}_f(x,y)) & = \begin{cases}
    \left[\widetilde{h}(\mc{S}_f(x,y))\right]^{1/\phi} & \text{if } \mc{S}_f(x,y)\in \widetilde{h}^{-1}\left([0,1]\right)\\
    \widetilde{h}(\mc{S}_f(x,y)) & \text{if } \mc{S}_f(x,y)\in \widetilde{h}^{-1}\left([1,\infty)\right)
\end{cases}\\
& = \begin{cases}
    \widetilde{h}(x,y) & \text{if } (x,y)\in \widetilde{h}^{-1}\left([0,1]\right)\\
    \left[\widetilde{h}(x,y)\right]^\phi & \text{if } (x,y)\in \widetilde{h}^{-1}\left([1,\infty)\right)
\end{cases} = \left[h_{f,z_0}(x,y)\right]^\phi,
\end{align*}
as desired.
\\

\noindent\boxed{\textbf{(2).}} Let $\|\cdot\| : \C^2\to \R$ be any norm on $\C^2$. Also let $\|\cdot\|_\infty : \C^2\to \R$ be the $l^\infty$-norm, namely
\[
\|(z,w)\|_\infty = \max\left\{|z|,|w|\right\}.
\]
By the equivalence of all norms on a given finite-dimensional $\C$-vector space, we know there exist positive constants $c,C > 0$ such that
\[
c\cdot \|(z,w)\|_\infty \leq \|(z,w)\| \leq C\cdot \|(z,w)\|_\infty
\]
for all $(z,w)\in \C^2$. We will use these constants to prove that 
\[
\lim_{n\to \infty} \left\|\mc{S}_f^{\circ n}(x,y) - (z_0,z_0)\right\|_\infty^{1/\phi^n} = \lim_{n\to \infty} \left\|\mc{S}_f^{\circ n}(x,y) - (z_0,z_0)\right\|^{1/\phi^n}
\]
for any $(x,y)\in \mc{A}(f,z_0)$ such that either limit exists. %(we do not yet know the result is $h_{f,z_0}(x,y)$). 
Before proceeding, note that since $\lim_{n\to\infty} 1/\phi^n = 0$, it follows that
\[
\lim_{n\to\infty} K^{1/\phi^n} = K^0 = 1
\]
for any positive constant $K > 0$. %Therefore, if $(x,y)\in \mc{A}(x,y)$ is such that $\lim_{n\to \infty} \left\|\mc{S}_f^{\circ n}(x,y) - (z_0,z_0)\right\|^{1/\phi^n}$ exists (including the case when $\|\cdot\| = \|\cdot\|_\infty$), we find
% \[
% \lim_{n\to\infty} \left[K\cdot \left\|\mc{S}_f^{\circ n}(x,y) - (z_0,z_0)\right\|\right]^{1/\phi^n} = \left(\lim_{n\to\infty} K^{1/\phi^n}\right)\cdot \left(\lim_{n\to \infty} \left\|\mc{S}_f^{\circ n}(x,y) - (z_0,z_0)\right\|^{1/\phi^n}\right) = \lim_{n\to \infty} \left\|\mc{S}_f^{\circ n}(x,y) - (z_0,z_0)\right\|^{1/\phi^n}
% \]

Suppose $(x_0,y_0)\in \mc{A}(f,z_0)$ is such that $\lim_{n\to \infty} \left\|\mc{S}_f^{\circ n}(x_0,y_0) - (z_0,z_0)\right\|_\infty^{1/\phi^n} = 0$. Then, for $K\in \{c,C\}$, we see
\[
\lim_{n\to\infty} \left[K\cdot \left\|\mc{S}_f^{\circ n}(x_0,y_0) - (z_0,z_0)\right\|_\infty\right]^{1/\phi^n} = \left(\lim_{n\to\infty} K^{1/\phi^n}\right)\cdot \left(\lim_{n\to \infty} \left\|\mc{S}_f^{\circ n}(x_0,y_0) - (z_0,z_0)\right\|_\infty^{1/\phi^n}\right) = 0.
\]
Since norm equivalence implies the inequality
\[
\left[c\cdot \left\|\mc{S}_f^{\circ n}(x_0,y_0) - (z_0,z_0)\right\|_\infty\right]^{1/\phi^n}\leq \left\|\mc{S}_f^{\circ n}(x_0,y_0) - (z_0,z_0)\right\|^{1/\phi^n} \leq \left[C\cdot \left\|\mc{S}_f^{\circ n}(x_0,y_0) - (z_0,z_0)\right\|_\infty\right]^{1/\phi^n}
\]
for all $n\in \Z^+$, the Squeeze Theorem applies to say $\lim_{n\to \infty} \left\|\mc{S}_f^{\circ n}(x,y) - (z_0,z_0)\right\|^{1/\phi^n}$ exists, and is equal to
\[
\lim_{n\to \infty} \left\|\mc{S}_f^{\circ n}(x,y) - (z_0,z_0)\right\|^{1/\phi^n} = 0 = \lim_{n\to \infty} \left\|\mc{S}_f^{\circ n}(x,y) - (z_0,z_0)\right\|_\infty^{1/\phi^n}
\]
in this case. 

Next, let $(x_1,y_1)\in \mc{A}(f,z_0)$ be such that $\lim_{n\to \infty} \left\|\mc{S}_f^{\circ n}(x_1,y_1) - (z_0,z_0)\right\|_\infty^{1/\phi^n} > 0$. Then it must be that $\left\|\mc{S}_f^{\circ n}(x_1,y_1) - (z_0,z_0)\right\|_\infty^{1/\phi^n} > 0$ for all $n\in \Z^+$, because otherwise
\begin{align*}
\left\|\mc{S}_f^{\circ n}(x_1,y_1) - (z_0,z_0)\right\|_\infty^{1/\phi^n} = 0 & \implies \left\|\mc{S}_f^{\circ n}(x_1,y_1) - (z_0,z_0)\right\|_\infty = 0 \\
& \implies \mc{S}_f^{\circ n}(x_1,y_1) = (z_0,z_0) \\
& \implies \mc{S}_f^{\circ m}(x_1,y_1) = \mc{S}_f^{\circ (m-n)}(z_0,z_0) = (z_0,z_0) \text{ for all } m\geq n\\
& \implies \left\|\mc{S}_f^{\circ m}(x_1,y_1) - (z_0,z_0)\right\|_\infty = 0 \text{ for all } m\geq n\\
& \implies \left\|\mc{S}_f^{\circ m}(x_1,y_1) - (z_0,z_0)\right\|_\infty^{1/\phi^m} = 0 \text{ for all } m\geq n\\
& \implies \lim_{k\to\infty} \left\|\mc{S}_f^{\circ k}(x_1,y_1) - (z_0,z_0)\right\|_\infty^{1/\phi^k} = 0,
\end{align*}
which would be a contradiction. Thus, we may always divide through by $\left\|\mc{S}_f^{\circ n}(x_1,y_1) - (z_0,z_0)\right\|_\infty^{1/\phi^n}$ in the norm equivalence inequality
\[
\left[c\cdot \left\|\mc{S}_f^{\circ n}(x_1,y_1) - (z_0,z_0)\right\|_\infty\right]^{1/\phi^n}\leq \left\|\mc{S}_f^{\circ n}(x_1,y_1) - (z_0,z_0)\right\|^{1/\phi^n} \leq \left[C\cdot \left\|\mc{S}_f^{\circ n}(x_1,y_1) - (z_0,z_0)\right\|_\infty\right]^{1/\phi^n}
\]
to obtain
\[
c^{1/\phi^n} \leq \frac{\left\|\mc{S}_f^{\circ n}(x_1,y_1) - (z_0,z_0)\right\|^{1/\phi^n}}{\left\|\mc{S}_f^{\circ n}(x_1,y_1) - (z_0,z_0)\right\|^{1/\phi^n}_\infty} \leq C^{1/\phi^n}.
\]
The Squeeze Theorem now applies to say
\[
\lim_{n\to\infty} \frac{\left\|\mc{S}_f^{\circ n}(x_1,y_1) - (z_0,z_0)\right\|^{1/\phi^n}}{\left\|\mc{S}_f^{\circ n}(x_1,y_1) - (z_0,z_0)\right\|^{1/\phi^n}_\infty} = 1,
\]
whence 
\begin{align*}
\lim_{n\to \infty} \left\|\mc{S}_f^{\circ n}(x,y) - (z_0,z_0)\right\|^{1/\phi^n} & = \left(\lim_{n\to\infty} \frac{\left\|\mc{S}_f^{\circ n}(x_1,y_1) - (z_0,z_0)\right\|^{1/\phi^n}}{\left\|\mc{S}_f^{\circ n}(x_1,y_1) - (z_0,z_0)\right\|^{1/\phi^n}_\infty}\right)\cdot\left(\lim_{n\to \infty} \left\|\mc{S}_f^{\circ n}(x,y) - (z_0,z_0)\right\|^{1/\phi^n}_\infty\right)\\
& = \lim_{n\to \infty} \left\|\mc{S}_f^{\circ n}(x,y) - (z_0,z_0)\right\|^{1/\phi^n}_\infty,
\end{align*}
as desired.

We did not use the actual formula for $\|\cdot\|_\infty$ anywhere in the above argument, only the norm equivalence inequality. Observe that the same inequality implies
\[
\frac1C\cdot \|(z,w)\| \leq \|(z,w)\|_\infty \leq \frac1c\cdot \|(z,w)\|
\]
for all $(z,w)\in \C^2$, which is another norm equivalence inequality, with $c$ and $C$ replaced with $1/C$ and $1/c$, respectively. Therefore, repeating the arguments of the previous two paragraphs, but with the roles of $\|\cdot\|$ and $\|\cdot\|_\infty$ swapped, we obtain that 
\[
\lim_{n\to \infty} \left\|\mc{S}_f^{\circ n}(x,y) - (z_0,z_0)\right\|^{1/\phi^n} \text{ exists } \implies \lim_{n\to \infty} \left\|\mc{S}_f^{\circ n}(x,y) - (z_0,z_0)\right\|^{1/\phi^n}_\infty \text{ exists}
\]
and in this situation,
\[
\lim_{n\to \infty} \left\|\mc{S}_f^{\circ n}(x,y) - (z_0,z_0)\right\|^{1/\phi^n} = \lim_{n\to \infty} \left\|\mc{S}_f^{\circ n}(x,y) - (z_0,z_0)\right\|^{1/\phi^n}_\infty,
\]
as desired.

So it remains to show that
\[
\lim_{n\to \infty} \left\|\mc{S}_f^{\circ n}(x,y) - (z_0,z_0)\right\|^{1/\phi^n}_\infty = h_{f,z_0}(x,y)
\]
for all $(x,y)\in \mc{A}(f,z_0)$, which we shall do by using the explicit formula for $\|\cdot\|_\infty$. Note that by definition of the temporary function $\widetilde{h}$ from (1), and the Böttcher-type modulus $\widehat{\mc{H}}_{f,z_0}$, we already know
\[
\widetilde{h}(x,y) = \widehat{\mc{H}}_{f,z_0}(x,y) \cdot |x-z_0|^{\phi/\sqrt5}\cdot |y-z_0|^{1/\sqrt5} = \lim_{n\to\infty} \left|(\pi_1\circ\mc{S}_f^{\circ n})(x,y) - (z_0,z_0)\right|^{1/\phi^n}.
\]
But it is evident %from our proposition concerning the iterates of $\mc{S}_f$ on $\mc{A}(f,z_0)$ 
Proposition 3 that 
\[
\left|(\pi_2\circ\mc{S}_f^{\circ (n+1)})(x,y) - z_0\right| = \left|(\pi_1\circ\mc{S}_f^{\circ n})(x,y) - z_0\right|
\]
for all $n\in\Z_{\geq0}$ and all $(x,y)\in \mc{A}(f,z_0)$. Therefore,
\begin{align*}
\lim_{n\to\infty} \left|(\pi_2\circ\mc{S}_f^{\circ n})(x,y) - z_0\right|^{1/\phi^n} & = \lim_{n\to\infty} \left|(\pi_2\circ\mc{S}_f^{\circ (n+1)})(x,y) - z_0\right|^{1/\phi^{n+1}} \\
& = \lim_{n\to\infty} \left|(\pi_1\circ\mc{S}_f^{\circ n})(x,y) - z_0\right|^{1/\phi^{n+1}}\\
& = \lim_{n\to\infty} \left[\left|(\pi_1\circ\mc{S}_f^{\circ n})(x,y) - z_0\right|^{1/\phi^n}\right]^{1/\phi}\\
& = \left[\lim_{n\to\infty} \left|(\pi_1\circ\mc{S}_f^{\circ n})(x,y) - z_0\right|^{1/\phi^n}\right]^{1/\phi}\\
& = \left[\widetilde{h}(x,y)\right]^{1/\phi}.
\end{align*}
We conclude
\begin{align*}
\lim_{n\to \infty} & \left\|\mc{S}_f^{\circ n}(x,y) - z_0\right\|^{1/\phi^n}_\infty \\
& = \lim_{n\to \infty} \left[\max\left\{\left|(\pi_1\circ\mc{S}_f^{\circ n})(x,y) - z_0\right|, \left|(\pi_2\circ\mc{S}_f^{\circ n})(x,y) - z_0\right|\right\}\right]^{1/\phi^n}\\
& = \lim_{n\to\infty} \max\left\{\left|(\pi_1\circ\mc{S}_f^{\circ n})(x,y) - z_0\right|^{1/\phi^n}, \left|(\pi_2\circ\mc{S}_f^{\circ n})(x,y) - z_0\right|^{1/\phi^n}\right\}\\
& = \max\left\{\lim_{n\to\infty} \left|(\pi_1\circ\mc{S}_f^{\circ n})(x,y) - z_0\right|^{1/\phi^n}, \lim_{n\to\infty} \left|(\pi_2\circ\mc{S}_f^{\circ n})(x,y) - z_0\right|^{1/\phi^n}\right\}\\
& = \max\left\{\widetilde{h}(x,y), \left[\widetilde{h}(x,y)\right]^{1/\phi}\right\}\\
& = h_{f,z_0}(x,y),
\end{align*}
as desired.
\\

\noindent\boxed{\textbf{(3).}} Certainly we have $\bigcup_{n=0}^\infty \mc{S}_f^{\circ(-n)}(\{(z_0,z_0)\}) \subseteq h_{f,z_0}^{-1}(\{0\})$: if $(x_0,y_0)\in \mc{A}(f,z_0)$ is such that $\mc{S}_f^{\circ N}(x_0,y_0) = (z_0,z_0)$ for some $N\in \Z_{\geq0}$, then 
\[
\mc{S}_f^{\circ m}(x_0,y_0) = \mc{S}_f^{\circ (m-N)}(z_0,z_0) = (z_0,z_0)
\]
for all $m\geq N$, hence
\[
h_{f,z_0}(x_0,y_0) = \lim_{m\to\infty} \left\|\mc{S}_f^{\circ m}(x_0,y_0)-(z_0,z_0)\right\|^{1/\phi^m} = \lim_{m\to\infty} \left\|(z_0,z_0)-(z_0,z_0)\right\|^{1/\phi^m} = 0,
\]
as desired.

Conversely, let $(x_1,y_1)\in h_{f,z_0}^{-1}(\{0\})$ be given. This implies $\widetilde{h}(x_1,y_1) = 0$, or in other words,
\[
\widehat{\mc{H}}_{f,z_0}(x_1,y_1) \cdot |x_1-z_0|^{\phi/\sqrt5} \cdot |y-z_0|^{1/\sqrt5} = 0.
\]
From here, we obtain three possible cases: $\widehat{\mc{H}}(x_1,y_1) = 0$, or $x_1 = z_0$, or $y_1 = z_0$. We show that all three cases imply $(x_1,y_1)\in \bigcup_{n=0}^\infty \mc{S}_f^{\circ(-n)}(\{(z_0,z_0)\})$.

The latter two cases can be dealt with using an explicit computation from the formula for $\mc{S}_f$ in Lemma \hyperlink{2:Jacobian}{2.1}: for all $x,y\in \C$ for which $(x,z_0)$ and $(z_0,y)$ are not indeterminate for $\mc{S}_f$, we find
\begin{align*}
\mc{S}_f(x,z_0) & = \begin{cases}
    \left(\dfrac{f(x)z_0 - f(z_0)x}{f(x)-f(z_0)}, x\right) & \text{if } f(x)\neq 0 \\[0.7em]
    \left(\dfrac{z_0f'(z_0)-f(z_0)}{f'(z_0)},z_0\right) & \text{if } x = z_0
\end{cases} = \begin{cases}
    (z_0,x) & \text{if } f(x)\neq 0\\
    (z_0,z_0) & \text{if } x = z_0
\end{cases}
\end{align*}
and
\begin{align*}
\mc{S}_f(z_0,y) & = \begin{cases}
    \left(\dfrac{f(z_0)y - f(y)z_0}{f(z_0)-f(y)}, z_0\right) & \text{if } f(y)\neq 0\\[0.7em]
    \left(\dfrac{z_0f'(z_0)-f(z_0)}{f'(z_0)},z_0\right) & \text{if } y = z_0
\end{cases} = (z_0,z_0).
\end{align*}
That is, $\mc{S}_f^{\circ 2}(x,z_0) = (z_0,z_0) = \mc{S}_f(z_0,y)$. So in our case, if $x_1 = z_0$ or $y_1 = z_0$, then $\mc{S}_f^{\circ 2}(x_1,y_1) = (z_0,z_0)$, as desired.

Suppose instead that $\widehat{\mc{H}}_{f,z_0}(x_1,y_1) = 0$. Recall from Theorem \hyperlink{3:BöttcherFunction}{3.1} that 
\[
\widehat{\mc{H}}_{f,z_0}(z_0,z_0) = \left|\mc{H}_{f,z_0,w_0}(z_0,z_0)\right| = \left|\frac{f''(z_0)}{2f'(z_0)}\right|^{\frac{5+3\sqrt5}{10}} > 0.
\]
Since $\lim_{n\to\infty} \mc{S}_f^{\circ n}(x_1,y_1) = (z_0,z_0)$ by definition of $\mc{A}(f,z_0)$, it follows from the continuity of $\widehat{\mc{H}}_{f,z_0}$ that there exists an $N\in \Z_{\geq0}$ such that
\[
n\geq N \text{ implies } \mc{S}_f^{\circ n}(x_1,y_1) \in \widehat{\mc{H}}_{f,z_0}^{-1}((0,\infty)).
\]
In particular, $N > 0$ since $\widehat{\mc{H}}_{f,z_0}(x_1,y_1) = 0$. Apply induction to the functional equation for $\widehat{\mc{H}}_{f,z_0}$ to now see
\[
0 = \widehat{\mc{H}}_{f,z_0}(x_1,y_1) = \left[\widehat{\mc{H}}_{f,z_0}\left(\mc{S}_f^{\circ N}(x_1,y_1)\right)\right]^{\phi^{-N}} \cdot \prod_{k=0}^{N-1}\left|\mc{G}_{f,z_0}\left(\mc{S}_f^{\circ k}(x_1,y_1)\right)\right|^{\phi^{-k}/\sqrt5},
\]
where as usual $\mc{G}_{f,z_0}$ is from Lemma \hyperlink{2:PowerSeries}{2.2}. Since $\widehat{\mc{H}}_{f,z_0}\left(\mc{S}_f^{\circ N}(x_1,y_1)\right) > 0$ by definition of $N$, it follows that there exists a $0\leq k\leq N-1$ such that $\mc{G}_{f,z_0}\left(\mc{S}_f^{\circ k}(x_1,y_1)\right) = 0$. Immediately this implies
\begin{align*}
\mc{S}_f^{\circ (k+1)}(x_1,y_1) & = \mc{S}_f\left(\mc{S}_f^{\circ k}(x_1,y_1)\right) \\
& = \left(z_0 + \mc{G}_{f,z_0}\left(\mc{S}_f^{\circ k}(x_1,y_1)\right)\left(\pi_1\left(\mc{S}_f^{\circ k}(x_1,y_1)\right)-z_0\right)\left(\pi_2\left(\mc{S}_f^{\circ k}(x_1,y_1)\right)-z_0\right), \pi_1\left(\mc{S}_f^{\circ k}(x_1,y_1)\right)\right)\\
& = \left(z_0, \pi_1\left(\mc{S}_f^{\circ k}(x_1,y_1)\right)\right),
\end{align*}
whence our work on the other cases implies
\[
\mc{S}_f^{\circ (k+2)}(x_1,y_1) = (z_0,z_0),
\]
as desired.

In all cases, we have proved that $(x_1,y_1)\in h_{f,z_0}^{-1}(\{0\})$ implies $(x_1,y_1)\in \bigcup_{n=0}^\infty \mc{S}_f^{\circ(-n)}(\{(z_0,z_0)\})$. We conclude
\[
h_{f,z_0}^{-1}(\{0\}) = \bigcup_{n=0}^\infty \mc{S}_f^{\circ(-n)}(\{(z_0,z_0)\}),
\]
as claimed.
\\

\noindent\boxed{\textbf{(4).}} Let $(x,y)\in \mc{A}(f,z_0)$ be given. Then $\lim_{n\to \infty} \mc{S}_f^{\circ n}(x,y) = (z_0,z_0)$ by definition of the basin. Since $h_{f,z_0}$ is continuous, this implies
\[
\lim_{n\to\infty} h_{f,z_0}\left(\mc{S}_f^{\circ n}(x,y)\right) = h_{f,z_0}\left(\lim_{n\to\infty} \mc{S}_f^{\circ n}(x,y)\right) = h_{f,z_0}(z_0,z_0) = 0
\]
(using (3) for the final equality). But on the other hand, induction on the functional equation for $h_{f,z_0}$ says
\[
h_{f,z_0}\left(\mc{S}_f^{\circ n}(x,y)\right) = \left[h_{f,z_0}(x,y)\right]^{\phi^n}.
\]
Because $\lim_{n\to\infty} \phi^n = \infty$, the only way we can have 
\[
\lim_{n\to\infty} \left[h_{f,z_0}(x,y)\right]^{\phi^n} = \lim_{n\to\infty} h_{f,z_0}\left(\mc{S}_f^{\circ n}(x,y)\right) = 0
\]
is if we have $0\leq h_{f,z_0}(x,y) < 1$, as desired.

Observe that for our temporary function $\widetilde{h}$, we now have
\[
0\leq \widetilde{h}(x,y) \leq \max\left\{\widetilde{h}(x,y), \left[\widetilde{h}(x,y)\right]^{1/\phi}\right\} = h_{f,z_0}(x,y) < 1.
\]
This implies $0\leq\widetilde{h}(x,y) < \left[\widetilde{h}(x,y)\right]^{1/\phi} < 1$ since $0<1/\phi < 1$. Therefore, the maximum is
\[
h_{f,z_0}(x,y) = \left[\widetilde{h}(x,y)\right]^{1/\phi} = \left[\widehat{\mc{H}}_{f,z_0}(x,y)\right]^{1/\phi}\cdot |x-z_0|^{1/\sqrt5}\cdot |y-z_0|^{\phi^{-1}/\sqrt5}
\]
for all $(x,y)\in \mc{A}(f,z_0)$, as claimed.
\end{proof}

The final form of $h_{f,z_0}$ given in Corollary \hyperlink{3:Potential}{3.2}(4) reveals that if we let $\log \circ h_{f,z_0}$ be the Green's function
\[
(\log \circ h_{f,z_0})(x,y) = \frac1\phi\cdot \ln\left(\widehat{\mc{H}}_{f,z_0}(x,y)\right) + \frac1{\sqrt5}\cdot \ln\left(|x-z_0|\right) + \frac{\phi^{-1}}{\sqrt5}\cdot \ln\left(|y-z_0|\right),
\]
then $\log \circ h_{f,z_0}$ is pluriharmonic in slit neighborhoods of $(z_0,z_0)$. More specifically, there exists a radius $r > 0$ such that $D_r(z_0) \times D_r(z_0) \subseteq \mc{A}(f,z_0)$, and for any pair of simple curves $\gamma_x, \gamma_y : [0,1] \to D_r(0)$ from $z_0$ to a point in the boundary circle, there exists a holomorphic map $\Gamma : D_r(0)\setminus \gamma_x([0,1]) \times D_r(0)\setminus \gamma_y([0,1]) \to \C$ such that
\[
\re(\Gamma(a,b)) = (\log \circ h_{f,z_0})(a+z_0,b+z_0)
\]
for all $(a,b)$ in the domain of $\Gamma$. This $\Gamma$ is obtained by choosing three branches of the complex logarithm: one defined in some disk about $\widehat{\mc{H}}_{f,z_0}(z_0,z_0)$ (to deal with $\ln \circ \widehat{\mc{H}}_{f,z_0}$),
and one each for the slit disks $D_r(0)\setminus \gamma_x([0,1])$ and $D_r(0)\setminus \gamma_y([0,1])$ (to deal with $\ln(|x-z_0|)$ and $\ln(|y-z_0|)$).

Using the functional equation for $h_{f,z_0}$, one can extend this argument to say $\log\circ h_{f,z_0}$ is in fact pluriharmonic (in the more general sense that it is locally the real part of a holomorphic map) at every point where it is finite, i.e., the set $h_{f,z_0}^{-1}(\{0\})$ of iterated preimages of $(z_0,z_0)$. 

Although $h_{f,z_0}$ is not real-analytic near $(z_0,z_0)$ itself, the Böttcher-type modulus $\widehat{\mc{H}}_{f,z_0}$ is (since $\widehat{\mc{H}}_{f,z_0}(z_0,z_0) \neq 0$). %More generally, each holomorphic germ $\mc{H}_{f,z_0,w_0}$ 
We can measure its derivatives by differentiating the holomorphic germs $\mc{H}_{f,z_0,w_0}$ at $(z_0,z_0)$. 

\begin{center}
\begin{minipage}{37em}
\noindent\textbf{Proposition \hyperlink{3:HPowerSeries}{3.3} (Derivative of $\mc{H}_{f,z_0,w_0}$)} \textit{Suppose $f''(z_0)\neq 0$, and take any logarithm $w_0\in \log\left(\frac{f''(z_0)}{2f'(z_0)}\right)$. If $\mc{H}_{f,z_0,w_0} : \mc{A}(f,z_0) \dashrightarrow \C$ is the associated holomorphic germ from Theorem \hyperlink{3:BöttcherFunction}{3.1}, then its complex Jacobian at $(z_0,z_0)$ is}
\[
D\mc{H}_{f,z_0,w_0}(z_0,z_0) = \frac1{\sqrt5}\cdot \frac{-3[f''(z_0)]^2 + 2f'(z_0)f'''(z_0)}{6f'(z_0)f''(z_0)} \cdot \mc{H}_{f,z_0,w_0}(z_0,z_0) \cdot \begin{bmatrix}
    \phi & 1
\end{bmatrix}.
\]
\end{minipage}
\end{center}
\begin{proof} 
Recall from the proof of item (1) in Theorem \hyperlink{3:BöttcherFunction}{3.1} that if $\mc{L} : D_{\frac12\cdot\left|\frac{f''(z_0)}{2f'(z_0)}\right|}\left(\frac{f''(z_0)}{2f'(z_0)}\right) \to \C$ is a holomorphic branch of the complex logarithm such that
\[
\mc{L}\left(\frac{f''(z_0)}{2f'(z_0)}\right) = w_0,
\]
then there exists an $r>0$ such that $\mc{S}_f(D_r(z_0)\times D_r(z_0)) \subseteq D_r(z_0)\times D_r(z_0)$, $\mc{G}_{f,z_0}$ is nonzero on $D_r(z_0)\times D_r(z_0)$, and $\mc{H}_{f,z_0,w_0}$ is given by
\[
\mc{H}_{f,z_0,w_0}(x,y) = \lim_{n\to\infty}  \exp\left(\sum_{k=0}^{n-1} \frac{F_{n-k}}{\phi^n}\cdot \mc{L}\left(\mc{G}_{f,z_0}(\mc{S}_f^{\circ k}(x,y))\right)\right)
\]
for all $(x,y)\in D_r(z_0)\times D_r(z_0)$. 

On this bi-disk, we can define a holomorphic branch of the map $(x,y)\mapsto \left[\mc{H}_{f,z_0,w_0}(\mc{S}_f(x,y))\right]^{1/\phi} \cdot \left[\mc{G}_{f,z_0}(x,y)\right]^{1/\sqrt5}$, namely
\begin{align*}
\left[\mc{H}_{f,z_0,w_0}(\mc{S}_f(x,y))\right]^{1/\phi} & \cdot \left[\mc{G}_{f,z_0}(x,y)\right]^{1/\sqrt5}\\
& = \lim_{n\to\infty} \exp\left(\frac1\phi\cdot\sum_{k=0}^{n-1} \frac{F_{n-k}}{\phi^n}\cdot \mc{L}\left(\mc{G}_{f,z_0}\left(\mc{S}_f^{\circ k}(\mc{S}_f(x,y))\right)\right)\right)\cdot \exp\left(\frac1{\sqrt5}\cdot \mc{L}(\mc{G}_{f,z_0}(x,y))\right).
\end{align*}
Direct substitution of $(z_0,z_0)$ into this branch reveals
\begin{align*}
\left[\mc{H}_{f,z_0,w_0}(\mc{S}_f(z_0,z_0))\right]^{1/\phi} & \cdot \left[\mc{G}_{f,z_0}(z_0,z_0)\right]^{1/\sqrt5}\\
& = \lim_{n\to\infty} \exp\left(\frac1\phi\cdot\sum_{k=0}^{n-1} \frac{F_{n-k}}{\phi^n}\cdot \mc{L}\left(\mc{G}_{f,z_0}\left(\mc{S}_f^{\circ k}(\mc{S}_f(z_0,z_0))\right)\right)\right)\cdot \exp\left(\frac1{\sqrt5}\cdot \mc{L}(\mc{G}_{f,z_0}(z_0,z_0))\right)\\
& = \lim_{n\to\infty} \exp\left(\frac1\phi\cdot\sum_{k=0}^{n-1} \frac{F_{n-k}}{\phi^n}\cdot \mc{L}(\mc{G}_{f,z_0}(z_0,z_0))\right)\cdot \exp\left(\frac1{\sqrt5}\cdot \mc{L}(\mc{G}_{f,z_0}(z_0,z_0))\right)\\
& = \lim_{n\to\infty} \exp\left(\frac1\phi\cdot\sum_{k=0}^{n-1} \frac{F_{n-k}}{\phi^n}\cdot \mc{L}\left(\frac{f''(z_0)}{2f'(z_0)}\right)\right)\cdot \exp\left(\frac1{\sqrt5}\cdot \mc{L}\left(\frac{f''(z_0)}{2f'(z_0)}\right)\right)\\
& = \lim_{n\to\infty} \exp\left(\frac1\phi\cdot\sum_{k=0}^{n-1} \frac{F_{n-k}}{\phi^n}\cdot w_0\right)\cdot \exp\left(\frac1{\sqrt5}\cdot w_0\right)\\
& = \exp\left(\left(\frac1\phi\cdot \left(\lim_{n\to\infty}\sum_{k=0}^{n-1}\frac{F_{n-k}}{\phi^n}\right)+\frac1{\sqrt5}\right)\cdot w_0\right) = \exp\left(\left(\frac1\phi\cdot \frac{5+3\sqrt5}{10} + \frac1{\sqrt5}\right)\cdot w_0\right)\\
& = \exp\left(\frac{5+3\sqrt5}{10}\cdot w_0\right) = \mc{H}_{f,z_0,w_0}(z_0,z_0).
\end{align*}
It also follows from items (3) and (4) of Theorem \hyperlink{3:BöttcherFunction}{3.1} that 
\[
\lim_{n\to \infty} \left|(\pi_1\circ \mc{S}_f^{\circ n})(x,y) - z_0\right|^{1/\phi^n} = \left|\left[\mc{H}_{f,z_0,w_0}(\mc{S}_f(z_0,z_0))\right]^{1/\phi} \cdot \left[\mc{G}_{f,z_0}(z_0,z_0)\right]^{1/\sqrt5}\right|\cdot |x-z_0|^{\phi/\sqrt{5}}\cdot |y-z_0|^{1/\sqrt5}
\]
for all $(x,y)\in D_r(z_0) \times D_r(z_0)$. Therefore, item (2) of Theorem \hyperlink{3:BöttcherFunction}{3.1} applies to say
\[
\left[\mc{H}_{f,z_0,w_0}(\mc{S}_f(x,y))\right]^{1/\phi} \cdot \left[\mc{G}_{f,z_0}(x,y)\right]^{1/\sqrt5} = \mc{H}_{f,z_0,w_0}(x,y)
\]
on the bi-disk $D_r(z_0) \times D_r(z_0)$.

We can differentiate the above functional equation, applying the chain and product rules to find the complex Jacobian is
\begin{align*}
D\mc{H}_{f,z_0,w_0}(x,y) & = D\left(\left[\mc{H}_{f,z_0,w_0}(\mc{S}_f(x,y))\right]^{1/\phi}\right)\cdot \left[\mc{G}_{f,z_0}(x,y)\right]^{1/\sqrt5} + \left[\mc{H}_{f,z_0,w_0}(\mc{S}_f(x,y))\right]^{1/\phi} \cdot D\left(\left[\mc{G}_{f,z_0}(x,y)\right]^{1/\sqrt5}\right)\\
& = \frac1\phi\cdot \left[\mc{H}_{f,z_0,w_0}(\mc{S}_f(x,y))\right]^{1/\phi-1}\cdot D\left(\mc{H}_{f,z_0,w_0}(\mc{S}_f(x,y))\right)\cdot \left[\mc{G}_{f,z_0}(x,y)\right]^{1/\sqrt5}\\
& \ \ \ \ \ + \left[\mc{H}_{f,z_0,w_0}(\mc{S}_f(x,y))\right]^{1/\phi}\cdot \frac1{\sqrt5}\cdot \left[\mc{G}_{f,z_0}(x,y)\right]^{1/\sqrt5-1}\cdot D\mc{G}_{f,z_0}(x,y)\\
& = \frac{\mc{H}_{f,z_0,w_0}(x,y)}{\phi\cdot \mc{H}_{f,z_0,w_0}(\mc{S}_f(x,y))} \cdot D\mc{H}_{f,z_0,w_0}(\mc{S}_f(x,y))D\mc{S}_f(x,y) + \frac{\mc{H}_{f,z_0,w_0}(x,y)}{\sqrt5\cdot \mc{G}_{f,z_0}(x,y)}\cdot D\mc{G}_{f,z_0}(x,y)
\end{align*}
on the bi-disk. Substituting $(x,y) = (z_0,z_0)$ reveals
\begin{align*}
D\mc{H}_{f,z_0,w_0}(z_0,z_0) & = \frac{\mc{H}_{f,z_0,w_0}(z_0,z_0)}{\phi\cdot \mc{H}_{f,z_0,w_0}(z_0,z_0)} \cdot D\mc{H}_{f,z_0,w_0}(z_0,z_0)D\mc{S}_f(z_0,z_0) + \frac{\mc{H}_{f,z_0,w_0}(z_0,z_0)}{\sqrt5\cdot \mc{G}_{f,z_0}(z_0,z_0)}\cdot D\mc{G}_{f,z_0}(z_0,z_0)\\
& = \frac1\phi\cdot D\mc{H}_{f,z_0,w_0}(z_0,z_0)\begin{bmatrix}
    0 & 0 \\
    1 & 0
\end{bmatrix} + \frac{\mc{H}_{f,z_0,w_0}(z_0,z_0)}{\sqrt5\cdot \mc{G}_{f,z_0}(z_0,z_0)}\cdot D\mc{G}_{f,z_0}(z_0,z_0)
\end{align*}
(where we applied Lemma \hyperlink{2:Jacobian}{2.1}). Rearranging, we obtain
\begin{align*}
\frac{\mc{H}_{f,z_0,w_0}(z_0,z_0)}{\sqrt5\cdot \mc{G}_{f,z_0}(z_0,z_0)}\cdot D\mc{G}_{f,z_0}(z_0,z_0) & = D\mc{H}_{f,z_0,w_0}(z_0,z_0) - \frac1\phi\cdot D\mc{H}_{f,z_0,w_0}(z_0,z_0)\begin{bmatrix}
    0 & 0 \\
    1 & 0
\end{bmatrix}\\
& = D\mc{H}_{f,z_0,w_0}(z_0,z_0)\left(\begin{bmatrix}
    1 & 0 \\
    0 & 1
\end{bmatrix}-\frac1\phi\cdot\begin{bmatrix}
    0 & 0 \\
    1 & 0
\end{bmatrix}\right) \\
& = D\mc{H}_{f,z_0,w_0}(z_0,z_0) \begin{bmatrix}
    1 & 0 \\
    -1/\phi & 1
\end{bmatrix} = D\mc{H}_{f,z_0,w_0}(z_0,z_0) \begin{bmatrix}
    1 & 0 \\
    1/\phi & 1
\end{bmatrix}^{-1}.
\end{align*}
Therefore,
\[
D\mc{H}_{f,z_0,w_0}(z_0,z_0) = \frac{\mc{H}_{f,z_0,w_0}(z_0,z_0)}{\sqrt5\cdot \mc{G}_{f,z_0}(z_0,z_0)}\cdot D\mc{G}_{f,z_0}(z_0,z_0)\begin{bmatrix}
    1 & 0 \\
    1/\phi & 1
\end{bmatrix}.
\]
It remains to compute $D\mc{G}_{f,z_0}(z_0,z_0)$.

Recall from the proof of Lemma \hyperlink{2:PowerSeries}{2.2} that 
\begin{align*}
\mc{G}_{f,z_0}(x,y) & = \mc{G}_{f,z_0}(z_0,z_0) + \frac1{(3-1)!\cdot 1!}\cdot \frac{\partial^3 \mc{S}_f^1}{\partial x^2\partial y}(z_0,z_0)\cdot (x-z_0) + \frac1{(3-2)!\cdot 2!}\cdot\frac{\partial^3\mc{S}_f^1}{\partial x\partial y^2}(z_0,z_0) \cdot (y-z_0) + \cdots\\
& = \mc{G}_{f,z_0}(z_0,z_0) + \frac12\cdot \frac{\partial^3 \mc{S}_f^1}{\partial x^2\partial y}(z_0,z_0)\cdot (x-z_0) + \frac12\cdot\frac{\partial^3\mc{S}_f^1}{\partial x\partial y^2}(z_0,z_0) \cdot (y-z_0) + \cdots,
\end{align*}
where $\mc{S}_f^1$ denotes the first coordinate function of $\mc{S}_f$.  Therefore, 
\[
D\mc{G}_{f,z_0}(z_0,z_0) = \frac12\cdot \begin{bmatrix}
    \dfrac{\partial^3 \mc{S}_f^1}{\partial x^2\partial y}(z_0,z_0) & \dfrac{\partial^3 \mc{S}_f^1}{\partial x\partial y^2}(z_0,z_0)
\end{bmatrix}.
\]
Also recall from Lemma \hyperlink{2:Jacobian}{2.1} that
\[
\frac{\partial \mc{S}_f^1}{\partial x}(z_0,y) = \begin{cases}
    \dfrac{f(y)-f'(z_0)(y-z_0)}{f(y)} & \text{if } y\neq z_0\\
    0 & \text{if } y = z_0
\end{cases} \text{ and } \frac{\partial \mc{S}_f^1}{\partial y}(x,z_0) = \begin{cases}
    \dfrac{f(x)-f'(z_0)(x-z_0)}{f(x)} & \text{if } x\neq z_0\\
    0 & \text{if } x = z_0
\end{cases}
\]
as long as $x$ and $y$ are sufficiently near $z_0$. So if $g : D_r(z_0) \dashrightarrow \C$ is the map
\[
g(z) = \begin{cases}
    \dfrac{f(z)-f'(z_0)(z-z_0)}{f(z)} & \text{if } z\neq z_0\\
    0 & \text{if } z = z_0
\end{cases} = \begin{cases}
    1-\dfrac{f'(z_0)(z-z_0)}{f(z)} & \text{if } z\neq z_0\\
    0 & \text{if } z = z_0
\end{cases},
\]
then $g$ is holomorphic near $z_0$ and
\[
\dfrac{\partial^3 \mc{S}_f^1}{\partial x^2\partial y}(z_0,z_0) = g''(z_0) =  \dfrac{\partial^3 \mc{S}_f^1}{\partial x\partial y^2}(z_0,z_0).
\]
Observe that
\[
g'(z) = \begin{cases}
    \dfrac{d}{dz}\left(1-\dfrac{f'(z_0)(z-z_0)}{f(z)}\right) & \text{if } z\neq z_0\\
    \dfrac{\partial^2\mc{S}_f}{\partial x\partial y}(z_0,z_0) & \text{if } z = z_0
\end{cases} = \begin{cases}
    -\dfrac{f'(z_0)f(z) - f'(z_0)(z-z_0)f'(z)}{[f(z)]^2} & \text{if } z\neq z_0\\
    \dfrac{f''(z_0)}{2f'(z_0)}  & \text{if } z = z_0
\end{cases}.
\]
So, the relevant derivative limit is 
\begin{align*}
g''(z_0) & = \lim_{z\to z_0} \frac{g'(z)-g'(z_0)}{z-z_0} = \lim_{z\to z_0} \frac1{z-z_0}\left[-\dfrac{f'(z_0)f(z) - f'(z_0)(z-z_0)f'(z)}{[f(z)]^2} - \frac{f''(z_0)}{2f'(z_0)}\right]\\
& = \lim_{z\to z_0} - \frac{2[f'(z_0)]^2\left(f(z) - f'(z)(z-z_0)\right) + f''(z_0)[f(z)]^2}{2f'(z_0)[f(z)]^2(z-z_0)}.
\end{align*}
Using the power series for $f$ about $z_0$, the numerator and denominator of the final fraction simplify to
\begin{align*}
2[f'(z_0)]^2&\left(f(z) - f'(z)(z-z_0)\right) + f''(z_0)[f(z)]^2\\
& = 2[f'(z_0)]^2\left(\left(f'(z_0)(z-z_0) + \frac12f''(z_0)(z-z_0)^2 + \frac16f'''(z_0)(z-z_0)^3 + O((z-z_0)^4)\right)\right.\\
& \ \ \ \ \ \ \ \ \ \ \ \ \ \ \ \ \ \ \ \ \ \ \ \ \ \ \ \left. - \left(f'(z_0) + f''(z_0)(z-z_0) + \frac12f'''(z_0)(z-z_0)^2 + O((z-z_0)^3)\right)(z-z_0)\right)\\
& \ \ \ \ + f''(z_0)\left[f'(z_0)(z-z_0) + \frac12f''(z_0)(z-z_0)^2 + O((z-z_0)^3)\right]^2\\
& = 2[f'(z_0)]^2\left(-\frac12f''(z_0)(z-z_0)^2 - \frac13f'''(z_0)(z-z_0)^3 + O((z-z_0)^4)\right)\\
& \ \ \ \ + f''(z_0)\left([f'(z_0)]^2(z-z_0)^2 + f'(z_0)f''(z_0)(z-z_0)^3 + O((z-z_0)^4)\right)\\
& = \left(-\frac23[f'(z_0)]^2f'''(z_0) + f'(z_0)[f''(z_0)]^2\right)(z-z_0)^3 + O((z-z_0)^4)
\end{align*}
and
\begin{align*}
2f'(z_0)[f(z)]^2(z-z_0) & = 2f'(z_0)\left[f'(z_0)(z-z_0) + O((z-z_0)^2)\right]^2(z-z_0) \\
& = 2[f'(z_0)]^3(z-z_0)^3 + O((z-z_0)^4),
\end{align*}
respectively. Factoring out a $(z-z_0)^3$, we conclude
\begin{align*}
g''(z_0) & = \lim_{z\to z_0} - \frac{2[f'(z_0)]^2\left(f(z) - f'(z)(z-z_0)\right) + f''(z_0)[f(z)]^2}{2f'(z_0)[f(z)]^2(z-z_0)}\\
& = \lim_{z\to z_0} - \frac{\left(-\frac23[f'(z_0)]^2f'''(z_0) + f'(z_0)[f''(z_0)]^2\right)(z-z_0)^3 + O((z-z_0)^4)}{2[f'(z_0)]^3(z-z_0)^3 + O((z-z_0)^4)}\\
& = \lim_{z\to z_0} - \frac{\left(-2[f'(z_0)]^2f'''(z_0) + 3f'(z_0)[f''(z_0)]^2\right) + O(z-z_0)}{6[f'(z_0)]^3 + O(z-z_0)}\\
& = \frac{2f'(z_0)f'''(z_0) - 3[f''(z_0)]^2}{6[f'(z_0)]^2}.
\end{align*}
This implies
\begin{align*}
D\mc{G}_{f,z_0}(z_0,z_0) & = \frac12\cdot \begin{bmatrix}
    \dfrac{\partial^3 \mc{S}_f^1}{\partial x^2\partial y}(z_0,z_0) & \dfrac{\partial^3 \mc{S}_f^1}{\partial x\partial y^2}(z_0,z_0)
    \end{bmatrix} = \frac{g''(z_0)}2\cdot\begin{bmatrix}
        1 & 1
    \end{bmatrix} \\
    & = \frac{- 3[f''(z_0)]^2 + 2f'(z_0)f'''(z_0) }{12[f'(z_0)]^2}\cdot\begin{bmatrix}
        1 & 1
    \end{bmatrix}.
\end{align*}
Finally, substituting into $D\mc{H}_{f,z_0,w_0}(z_0,z_0)$, we obtain
\begin{align*}
D\mc{H}_{f,z_0,w_0}(z_0,z_0) & = \frac{\mc{H}_{f,z_0,w_0}(z_0,z_0)}{\sqrt5\cdot \mc{G}_{f,z_0}(z_0,z_0)}\cdot D\mc{G}_{f,z_0}(z_0,z_0)\begin{bmatrix}
    1 & 0 \\
    1/\phi & 1
\end{bmatrix} \\
& = \frac{\mc{H}_{f,z_0,w_0}(z_0,z_0)}{\sqrt5\cdot \mc{G}_{f,z_0}(z_0,z_0)}\cdot \frac{- 3[f''(z_0)]^2 + 2f'(z_0)f'''(z_0) }{12[f'(z_0)]^2}\cdot\begin{bmatrix}
        1 & 1
    \end{bmatrix}\begin{bmatrix}
    1 & 0 \\
    1/\phi & 1
\end{bmatrix}\\
& = \frac{\mc{H}_{f,z_0,w_0}(z_0,z_0)}{\sqrt5}\cdot \frac{2f'(z_0)}{f''(z_0)}\cdot \frac{- 3[f''(z_0)]^2 + 2f'(z_0)f'''(z_0) }{12[f'(z_0)]^2}\cdot\begin{bmatrix}
    1+\dfrac1\phi & 1
\end{bmatrix}\\
& = \frac1{\sqrt5}\cdot \frac{- 3[f''(z_0)]^2 + 2f'(z_0)f'''(z_0) }{6f'(z_0)f''(z_0)}\cdot \mc{H}_{f,z_0,w_0}(z_0,z_0) \cdot \begin{bmatrix}
    \phi & 1
\end{bmatrix},
\end{align*}
as desired.
\end{proof}

%%%%
\section{Computer calculation of $h_{f,z_0}$}
\begin{figure}[h]%
\centering
\includegraphics[width=1.0\textwidth]{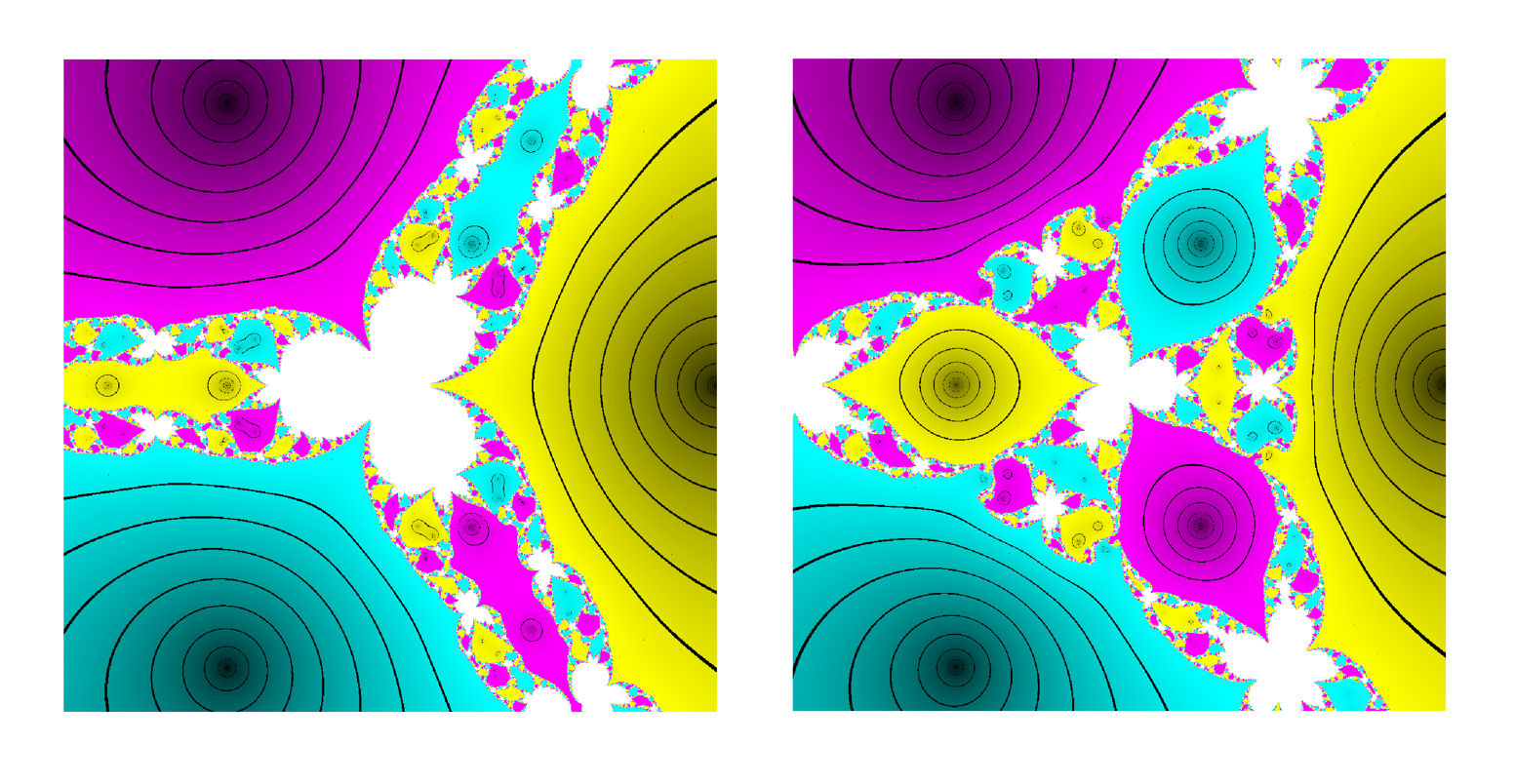}
\caption{Potential (hues) and contours (black) for the secant map of $f(z) = z^3-1$ in two planes. The white pixels do not converge to a fixed point. However, the center white regions correspond to the basin of an attracting 3-cycle discussed in \cite{BedfordFriggeSecant}.  (Left) Viewing the complex line $x = y$ parameterized by $x\in [-1,1] + i\cdot[-1,1]$. (Right) Viewing the critical curve $x+2y = 0$, parameterized by $y\in [-1,1]+ i\cdot [-1,1]$.}\label{fig1}
\end{figure}

\begin{figure}[h!]%
\centering
\includegraphics[width=1.0\textwidth]{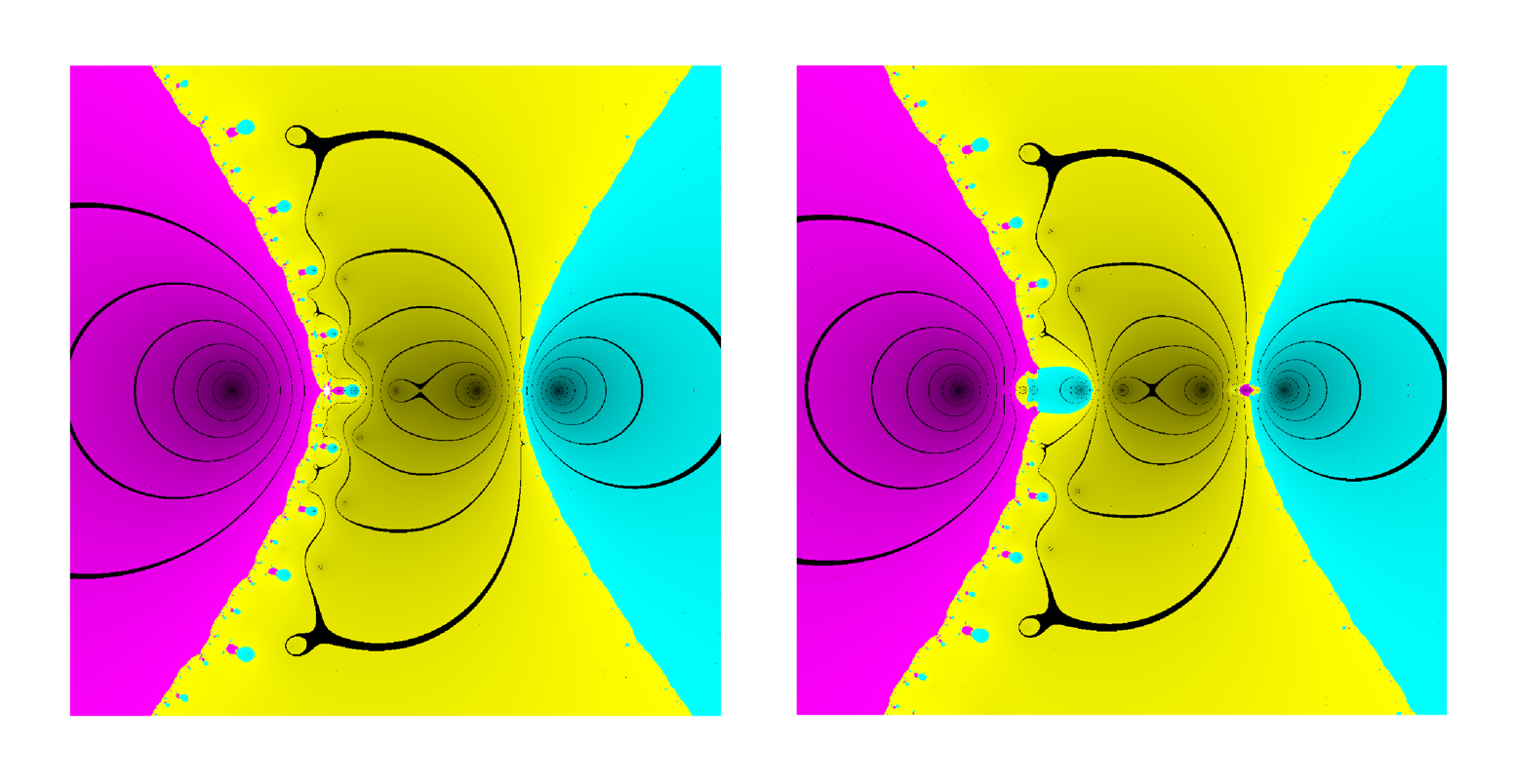}
\caption{Potential and contours for the secant map of $f(z) = (z^2-1)(z-1/2)$ in two planes.  (Left) Viewing the complex line $x = y$, parameterized by $x\in [-2,2]+i\cdot [-2,2]$. (Right) Viewing the critical curve $x+2y = 1/2$, parameterized by $y\in[-2,2]+i\cdot [-2,2]$.}\label{fig2}
\end{figure}

\begin{figure}[h!]%
\centering
\includegraphics[width=1.0\textwidth]{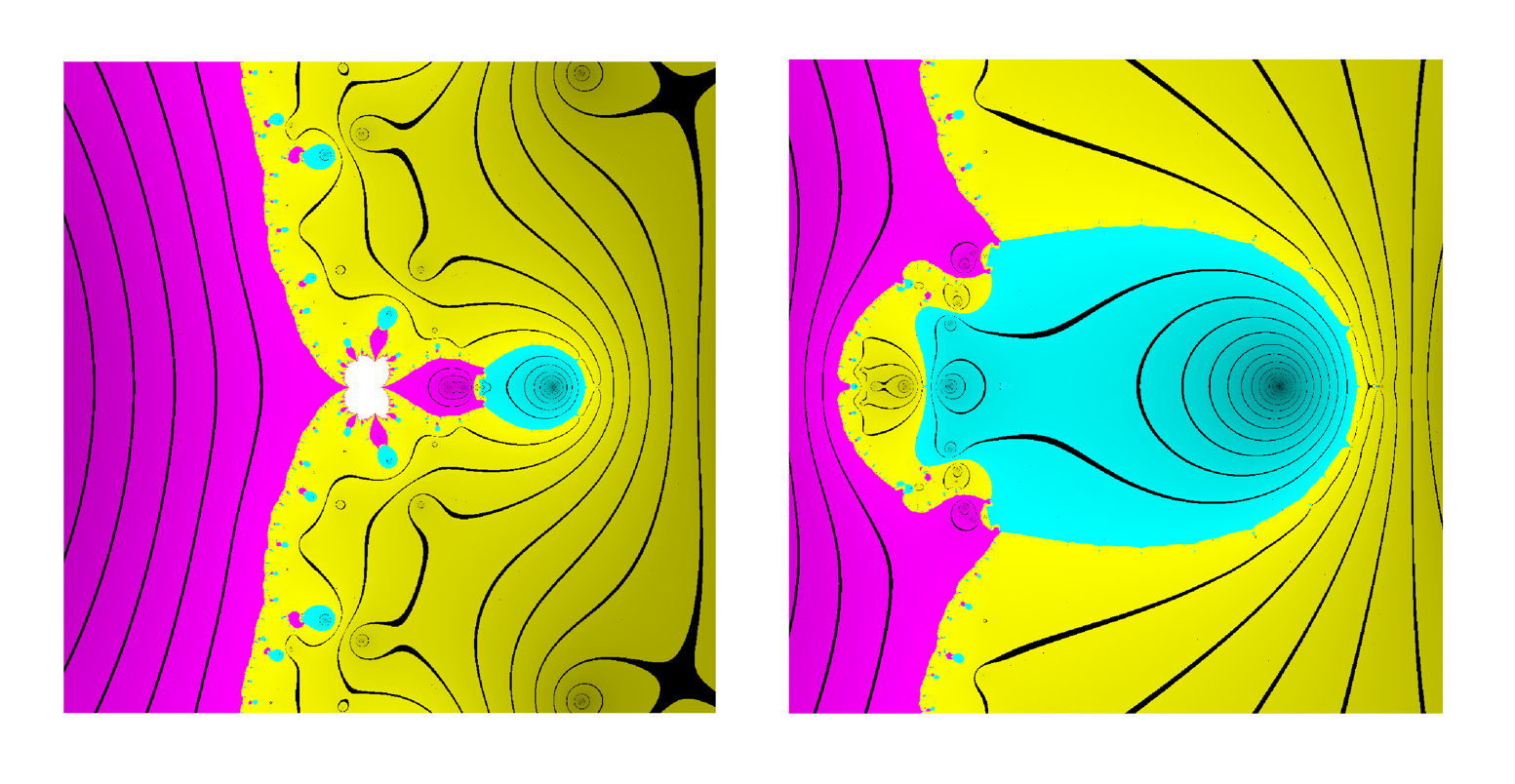}
\caption{Zoom in of both plots in Figure \ref{fig2}, with parameter ranging in $[-0.7,-0.1]+i\cdot[-0.3,0.3]$.}\label{fig3}
\end{figure}

Given $f$ and its set of simple zeros, we can use a computer to iterate $\mc{S}_f$ and detect within which attracting basin a point $(x,y)$ lies. With this done, we can calculate the value of the corresponding potential function. %There are some challenges with visualizing these data in a useful way. First and foremost is the problem of dimension: a complex $2$-dimensional domain such as that of $\mc{S}_f$ requires four real parameters to describe completely.
To visualize these data on a (real) two-dimensional image, we slice the domain $U\times U$ with various complex curves (typically complex lines) and color the points on the curve with hue according to their point of convergence (if it exists), and brightness according to the potential (since $0\leq h_{f,z_0}(x,y) < 1$). This is done in Figures \ref{fig1} through \ref{fig3} for a pair of cubics. 

There are always at least two particular curves of interest to plot. One is the diagonal plane $x = y$, which always contains the fixed points. Another is the critical curve
\[
\frac{f(x)-f(y)}{x-y} = f'(y)
\]
(minus its irreducible component $x = y$). This curve is the locus of points where the secant line to the graph of $f$ from $(x,f(x))$ to $(y,f(y))$ is in fact a tangent line, and tangent at $(y,f(y))$. For example, if $f$ is the cubic polynomial $f(z) = (z-a)(z-b)(z-c)$, then the associated critical curve is the line $x+2y = a+b+c$. 

\begin{figure}[h]%
\centering
\includegraphics[width=1.0\textwidth]{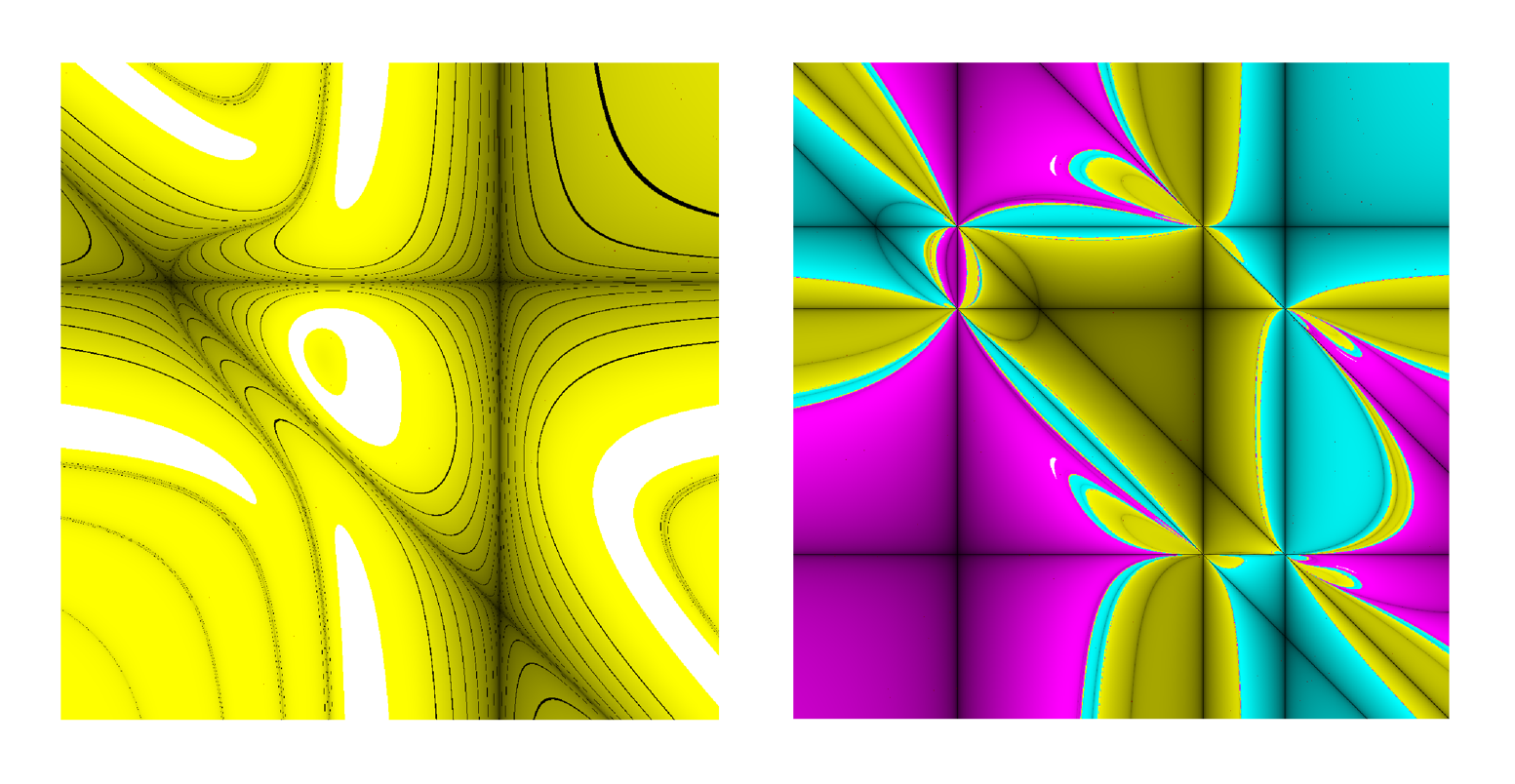}
\caption{Real dynamical plane for the secant maps of $f(z) = z^3-1$ (left) and $f(z) = (z^2-1)(z-1/2)$ (right), with $(x,y)$ ranging in $[-3,3]\times[-3,3]$ and $[-2,2]\times [-2,2]$, respectively. In the left figure, only one basin color appears--yellow, corresponding to the root $z_0 = 1$--because $\R^2$ is invariant under $\mc{S}_f$. In both figures, there are black $+$ shapes centered at the real fixed points. This is because the $\mc{S}_f$ collapses those vertical lines onto the fixed point, and takes the horizontal lines onto the vertical. Also in both figures, the white regions correspond to points in the basin of a 3-cycle; cf. \cite{GarijoJarqueRealSecant}.}\label{fig4}
\end{figure}

If we have $f(U\cap \R) \subseteq \R$ (i.e., if $f$ is the holomorphic extension of a smooth map $\R\dashrightarrow\R$), then one more avenue of visualization is available to us: the secant map $\mc{S}_f$ becomes the extension of a dynamical system $\R^2 \dashrightarrow \R^2$, and we can visualize the dynamical (real) plane of this underlying restriction. Note, however, that since $\R^2 \cap (U\times U)$ is invariant under $\mc{S}_f$ in this case, we will only be able to graph the potential for fixed points at diagonal real roots. This phenomenon is depicted in Figure \ref{fig4}.

% When using a computer to calculate the potential $h_{f,z_0}(x,y)$ (where $z_0$ is known beforehand), a heuristic algorithm immediately has access to the iteration number $N\in \Z_{\geq0}$ and the iterate $\mc{S}_f^{\circ N}(x,y)$. But we will prefer not to compute $h_{f,z_0}(x,y)$ via the approximation $\|\mc{S}_f^{\circ N}(x,y) - (z_0,z_0)\|^{1/\phi^N}$ due to possible catastrophic cancellation in the subtraction $\mc{S}_f^{\circ n}(x,y)-(z_0,z_0)$ (that is, two close numbers may have a difference significantly smaller than the machine's precision). %The problem is further compounded when we raise to the very small power $1/\phi^N$. 
% % because if $N$ is very large, then the approximation becomes as the indeterminate form $0^0$, which can be expensive to resolve accurately.

%%%%
\bibliographystyle{unsrt}  
\bibliography{references}  

\begin{thebibliography}{1}

\bibitem{DiezSecantConvergence}
P.~D\'iez.
\newblock A note on the convergence of the secant method for simple and multiple roots.
\newblock {\em Appl. Math. Lett.}, 16(8):1211--1215, 2003.

\bibitem{BedfordFriggeSecant}
Eric Bedford and Paul Frigge.
\newblock The secant method for root finding, viewed as a dynamical system.
\newblock {\em Dolomites Res. Notes Approx.}, 11(4):122--129, 2018.

\bibitem{HubbardPapadopolPotential}
John~H. Hubbard and Peter Papadopol.
\newblock Superattractive fixed points in {${\bf C}^n$}.
\newblock {\em Indiana Univ. Math. J.}, 43(1):321--365, 1994.

\bibitem{BuffEpsteinKochBottcherCoords}
Xavier Buff, Adam~L. Epstein, and Sarah Koch.
\newblock B\"ottcher coordinates.
\newblock {\em Indiana Univ. Math. J.}, 61(5):1765--1799, 2012.

\bibitem{GarijoJarqueRealSecant}
Antonio Garijo and Xavier Jarque.
\newblock Global dynamics of the real secant method.
\newblock {\em Nonlinearity}, 32(11):4557--4578, 2019.

\end{thebibliography}

\end{document}